\documentclass[a4paper,12pt,twoside]{article}
\usepackage[english]{babel}
\usepackage{amsmath}
\usepackage{amssymb}
\usepackage{latexsym}
\usepackage[mathscr]{euscript}
\usepackage{verbatim}
\usepackage{graphicx}
\DeclareGraphicsExtensions{.png}
\pdfoutput=1

\begin{document}
\vspace*{14mm}
\thispagestyle{empty} \large \noindent
\textbf{Generators of projective \vspace{6mm} MV-algebras} \\
\normalsize
\textbf{Francesco Lacava $\;\cdot\;$ Donato \vspace{50mm} Saeli} \\
\small
\textsc{Abstract } In the last decade, interest in projective MV-algebras has grown greatly; 
see [1], [5] e [6]. In this paper we establish a necessary and sufficient condition for $n$ elements 
of the free $n$-generator MV-algebra to generate a projective MV-algebra. This generalizes the 
characterization of the $n$ free generators proved in [7]. Using this, some classes of projective 
generators for bigenerated MV-algebras, are given. In particular, some effective procedures to 
determine, by elementary methods, generators of projective MV-algebras are \vspace{3mm} explained. \\
Nell'ultimo decennio, l'interesse per le MV-algebre proiettive \`e cresciuto 
sensibilmente; cfr. [1], [5] e [6]. In questa nota viene dimostrata una condizione necessaria e 
sufficiente perch\`e $n$ elementi di $\, Free_n. \,$ siano generatori di MV-algebre proiettive, che 
generalizza la caratterizzazione degli $n$ generatori liberi dimostrata in [7]. Utilizzando tale 
condizione nel caso di MV-algebre bigenerate sono state determinate alcune classi di generatori 
proiettivi. In particolare vengono esposti alcuni procedimenti effettivi per la determinazione dei 
generatori di algebre proiettive, con metodi \vspace{3mm} elementari. \\
\textsc{Keywords } MV-algebras $\;\cdot\;$ projective $\;\cdot\;$ retract $\;\cdot\;$ piecewise linear \vspace{3mm} map \\
\textsc{MSC } \vspace{6mm} 06D35 $\;\cdot\;$ 08B30 \footnote{\hspace*{-7.5mm} Francesco Lacava \\
Universit\`a degli Studi di Firenze \\
Dipartimento di Matematica ``Ulisse Dini'' \\
Viale Morgagni, 67/a - 50134 Firenze Italy \\
Tel.: +39-055-7877316 \\
Email: \vspace{2mm} francesco.lacava@math.unifi.it \\
Donato Saeli \\
Via Giovanni XXIII, 29 - 85100 Potenza, Italy \\
Tel.: +39-0971-51280 \\
Email: donato.saeli@gmail.com}
\newpage \pagestyle{myheadings} \markboth{\textsc{francesco lacava, donato saeli}}{\textsc{generators of projective mv-algebras}} \normalsize \noindent
\textbf{1 \vspace{2mm} Introduction} \\
Throughout this paper, by an MV-algebra $A$ we mean a \textit{semisimple} \\
MV-algebra, that is, a subalgebra of a direct power of the standard \\
MV-algebra [0, 1]. Thus for some set $T,$ without loss of generality we can assume that each $a \in 
A$ is a [0, 1]-valued function defined on $T.$

In particular, since the free $n$-generator MV-algebra $Free_n$ is semisimple, throughout we will 
identify $Free_n$ with the algebra of McNaughton functions on the $n$-cube [0, 1]$^n.$

An MV-algebra $A$ is said to be \textit{essentially} $n$-generated if it has a generating set of 
$n$ elements, but no generating set of $n$ - 1 elements.
We refer to [4], for all unexplained notation and terminology.

Throughout this paper, whenever $a$ generates $A$ we will assume that \\
$a(0) = 0.$

Let \textit{\L}$_{n+1}$ denote the {\L}ukasiewicz MV-chain with n+1 elements. Then following
[7], we will introduce the partial map $t:\,$\textit{\L}$_{n+1}\rightarrow\,$\textit{\L}$_{n+1}$ by 
stipulating that, for each $a \in$ \textit{\L}$_{n+1}$
$$t(a)=\left\{
\begin{array}{ll}
(ra)^{\prime } & \text{if there is} \ \; r \in \mathbb{Z}, \ \; r>0, \\
& \text{such that} \ ra<1 \ \; \text{and} \ \; \vspace{6mm} (r+1)a=1, \\
\text{undefined} & \text{otherwise.}
\end{array}
\right. $$
In particular, the map $t$ is undefined for $a=0, \, 1$. \ We further \vspace{.5mm} define \\
\hspace*{14mm} $t^0(a)=a$, \vspace{.5mm} \\
\hspace*{11mm} $t^{s+1}(a)=t(t^s(a))$. \vspace{1mm} \\
\textsc{Definition.} - An element  $a \in$ \L$_{n+1}$ is said to be a \textit{cyclic generator} of \L$_{n+1}, \,$ \\
if there is an integer $k \geq 0 \,$ such that \vspace{3mm} $\, t^k(a) = 1 / n.$ \\
\textsc{Proposition 1.1.} - \it If p is a prime number, then every element \ a$\, \in \, $\L$_{p+1}, \\
a \neq 0, \, 1, \ $ is a cyclic generator of \L$_{p+1}.$ \rm \hspace{\stretch{1}} This is [7, 2.4].\vspace{3mm}

For $p$ a prime number, let the integers $m$ and $p$ be such that $\,0 < m < p. \,$ Let $k$ be the smallest positive integer satisfying \ $t^k\Big(\dfrac{m}{p}\Big)=\dfrac{1}{p}$. \\
Let the sequence of elements of \textit{\L}$_{p+1},$ be given by
$$t\Big(\dfrac mp\Big)=\Big(n_1\dfrac mp\Big)' \! , \ t^2\Big(\dfrac mp\Big)=\Big(n_2t\Big(\dfrac mp\Big)\Big)' \! , \,\dots \, , \,t^k \Big(\dfrac mp\Big)=\Big(n_kt^{k-1}\Big(\dfrac mp\Big)\Big)'$$
Then the term $\, \gamma_{m,p} \,$ in the variable $x$ is defined by
$$\gamma_{m,p}(x) = \Big( n_k \big( n_{k-1} \, \dots \, \big( n_2 (n_1 x )' \big)' \, \dots \, \big)' \Big)'.$$
When interpreted in the $\, Free_1 \,$ algebra, the term $\, \gamma_{m,p} \,$ represents \\
a McNaughton function $\, g_{m,p} \,$ whose graph has three linear pieces, \\
two of which (eventually degenerate into points) are horizontal, respectively at level 0 and at level 1. Moreover, $\, g_{m,p}(m/p) = 1/p \,$ and $\, g_{m,p}(x) \neq 1/p \,$ for each $\, x \in [0,1] \,$ different from $\, m/p.$ \hspace{\stretch{1}} This is [7, 2.5].
\begin{center}
\includegraphics[scale =.3]{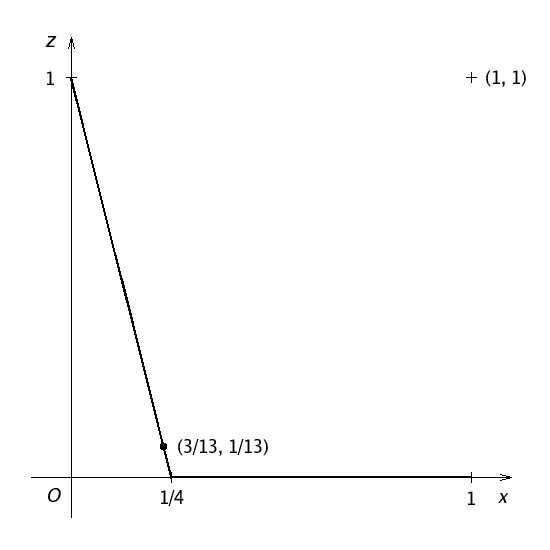} \hspace{22mm} \includegraphics[scale =.3]{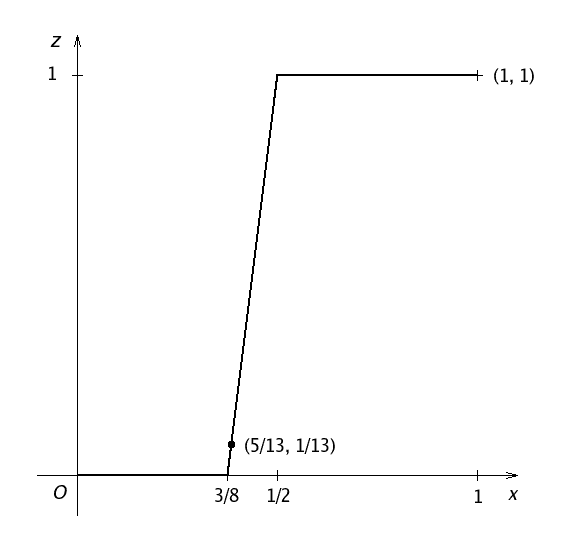} \\
$\, g_{_{3,13}} \,$ \hspace{60mm} \vspace{1mm} $\, g_{_{5,13}} \,$ \\
Fig. 1
\end{center}
For $p,$ a prime number, we can define the term $\, \lambda_p \,$ as follows:
$$\lambda_p(x) = \Big(px \wedge p \big((p-1)x \big)'\Big)'$$
If interpreted in the $\, Free_1 \,$ algebra, the term $\, \lambda_p \,$ represents a McNaughton function $\, l_p \,$ such that $\, l_p(x) = 0 \,$ if and only if $\, x = 1/p. \,$ The last nontrivial statement follows by noting that $\, x > 1/p \,$ implies \vspace{4mm} $\, p[1-(p-1)x] < 1.$
\begin{center}
\includegraphics[scale =.34]{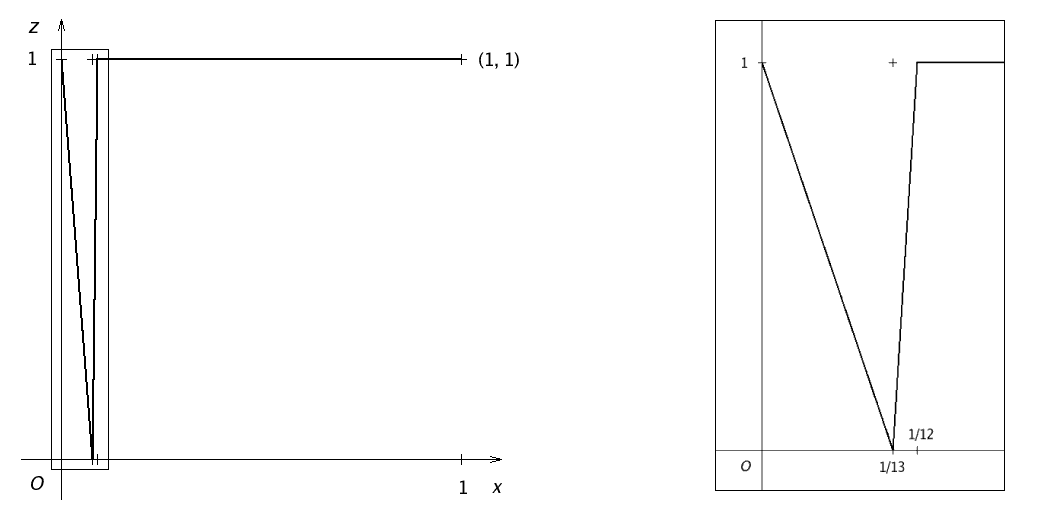} \\
$l_{_{13}}$ \vspace{2mm}\\
Fig. 2
\end{center}
For all integers $\,0 < m < p, \,$ with $p$ a prime number, let
$$\, \eta_{m,p}(x) = \lambda_p(\gamma_{m,p}(x)).$$
If interpreted in the $\, Free_1 \,$ algebra, the term $\, \eta_{m,p} \,$ represents a McNaughton function $\, t_{m,p} \,$ such \vspace{1mm} that: \\
(1.1) \, \it For all $\, x \in [0,1], \ \; t_{m,p}(x) = 0 \,$ if and only if \vspace{.5mm} $\, x = m/p.$ \rm \\
(1.2) \, \it For every $\, L \in \mathbb{Z}, \ \; L > 0 \,$ there is $\, Q = Q(L) \in \mathbb{Z}, \ \; Q > 0 \,$ such \vspace{.5mm} that, \\
\hspace*{10mm} for all $\, x \in [0,1] \,$ with $\, |x-m/p| \geq 1/L, \,$ we \vspace{.5mm} have $\, t_{m,p}(x) \geq 1/Q.$ \\
\rm \hspace*{\stretch{1}} This is \vspace{1mm} [7, 2.6].
\begin{center}
\includegraphics[scale =.44]{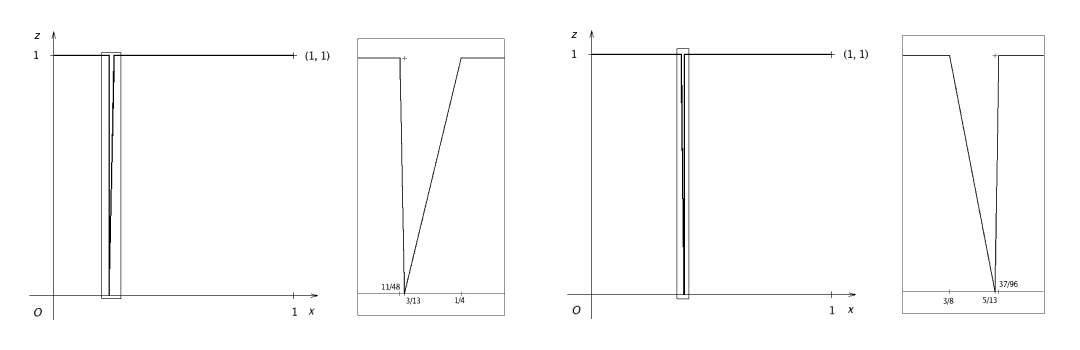} \\
$\, t_{_{3,13}} \,$ \hspace{60mm} \vspace{1mm} $\, t_{_{5,13}} \,$ \\
Fig. 3 
\end{center}
Following [4, 6.2.3], we say that an element $\, a \,$ of an MV-algebra \\
is archimedean if there is an integer $\, n \geq 1 \,$ such that $\, na + na = na.$ \vspace{2mm} \\
\textsc{Lemma 1.2.} - \it Suppose the elements $\, a_1, \, a_2, \dots \, , a_n \,$ $($essentially$)$ generate \\
the MV-algebra A  of $[0, 1]$-valued functions . Let $\, m_i < p_i \,$ be integers $\, > 0 \,$ with each $\, p_i \,$ prime $\, ( 1 \leq i \leq n ). \,$ Then the element
$$\bigvee_{i=1}^n \eta_{m_i,p_i}(a_i)$$
fails to be archimedean if and only if for every $\, \varepsilon >0 \,$ there is $\, \bar{\textit{\textbf{x}}} \in [0,1]^n \,$ \vspace{1.5mm} such that for $\, i=1, \dots, n, \,$ we \vspace{1.5mm} have $\, \Big|a_i(\bar{\textit{\textbf{x}}}) - \dfrac{m_i}{p_i} \Big| < \varepsilon$. \rm \\
\textsc{Proof.} - This follows at once from point (1.1) aforesaid. \vspace{2mm} \hspace{\stretch{1}} $\blacktriangle$ \\
\textsc{Proposition 1.3.} - \it Suppose the MV-algebras $\, A \ and \ B \,$ of $[0, \, 1]$-valued functions are essentially generated by sets $\, \{a_1,a_2,...,a_n\} \,$ and $\, \{b_1,b_2,...,b_n\} \,$ respectively. Then the following two conditions are equivalent: \\
\hspace*{2mm} $\bullet \ $ The map $\, f \, : \, a_i \mapsto b_i \,$ is $($uniquely$)$ extendible to an isomorphism of A \\ \hspace*{7mm} onto B \\
\hspace*{2mm} $\bullet \ $ For all integers $\, 0<m_i<p_i, \,$ with each $\, p_i \,$ prime, the element \\
\hspace*{4mm} $\ \displaystyle \bigvee_{i=1}^n \eta_{m_i,p_i}(a_i) \,$ is archimedean if and only if so is the element $\ \displaystyle \bigvee_{i=1}^n \eta_{m_i,p_i}(b_i).$ \rm \vspace{2mm} \\
\textsc{Proof.} - ($\Rightarrow $) The statement is obvious. \\
($\Leftarrow $) Vice versa, by contradiction, if we suppose that $\, f \,$ cannot be extended \\
to an isomorphism between $\, A $ and $ B, \,$ then a term $\, t \,$ exists so \vspace{.5mm} that \\
\hspace*{25mm} $\, t(a_1, a_2, \dots , a_n)= \hat{0} \,$ and \vspace{.5mm} $\, t(b_1, b_2, \dots , b_n) \neq \hat{0}. \,$ \\
It follows that there are $\, \textbf{\textit{x}}_0 \in [0,1]^n \,$ and $\, \delta >0 \,$ such that for every \\
$\, \textbf{\textit{x}} \in \pmb{\mathscr{I}}(\textbf{\textit{x}}_0, \delta) \,$ (the $\delta$-neighbourhood of $\, \textbf{\textit{x}}_0$) we \vspace{.5mm} have \\
\hspace*{30mm}$\, t(b_1(\textbf{\textit{x}}), b_2(\textbf{\textit{x}}), \dots , b_n(\textbf{\textit{x}})) \neq 0.$ \vspace{.5mm} \\
Since the $\, b_i(\textbf{\textit{x}}) \,$ are McNaughton functions, it is possible to choose a $\, \textbf{\textit{x}}_0 \,$ such that a neighborhood of $\, (b_1(\textbf{\textit{x}}_0), b_2( \textbf{\textit{x}}_0), \dots , b_n(\textbf{\textit{x}}_0)) \,$ is contained in the image by $\, \textbf{ \textit{b}}(\textbf{\textit{x}}) =(b_1(\textbf{\textit{x}}), b_2(\textbf{\textit{x}}), \dots , b_n(\textbf{\textit{x}})) \,$ of \vspace{.5mm} $\, \pmb{\mathscr{I}}(\textbf{\textit{x}}_0, \delta).$  \\
Furthermore there are $\, \bar{\textbf{\textit{x}}} \in [0,1]^n \cap \mathbb{Q} \,$ and $\, \bar{\delta} >0 \,$ such that for \vspace{1mm} every \\
$\, \textbf{\textit{x}} \in \pmb{\mathscr{I}}(\bar{\textbf{\textit{x}}}, \bar{\delta}) \,$ we have $\, t(b_1(\textbf{\textit{x}}), b_2(\textbf{\textit{x}}), \dots , b_n(\textbf{\textit{x}})) \neq 0, \ $ also $\, b_i(\bar{\textbf{\textit{x}}})=\dfrac{m_i}{p_i} \,$ with $\, p_i \,$ prime numbers. \ So clearly, the element $\, \displaystyle\bigvee_{i=1}^n \eta_{m_i,p_i}(b_i) \,$ is not archimedean and, by hypothesis, $\, \displaystyle\bigvee_{i=1}^n \eta_{m_i,p_i}(a_i) \,$ is not archimedean \vspace{1mm} either. \\
Now by the lemma 1.2 there is a $\, \bar{\varepsilon} \,$ such that, for each \textbf{\textit{y}} \vspace{1mm} satisfying \\
$\, \Big| y_i-\dfrac{m_i}{p_i} \Big| < \bar{\varepsilon} \ \; (i=1, \, \dots, \, n), \ $ we have $\, t(\textbf{\textit{y}}) \neq 0.$ \\
Since there is also an $\, \textbf{\textit{x}} \in [0,1]^n \,$ with $\, \Big| a(x_i)-\dfrac{m_i}{p_i} \Big| < \bar{\varepsilon} \ \; (i=1, \, \dots, \, n), \,$ \vspace{1mm} \\
we conclude $\, t(a_1(\textbf{\textit{x}}), a_2(\textbf{\textit{x}}), \dots , a_n(\textbf{\textit{x}})) \neq 0, \,$ which is a \vspace{1.5mm} contradiction. \hspace{\stretch{1}} $\blacktriangle$ \\
\indent In particular, the previous proposition tells us that in the in the free one-generator MV-algebra $\, Free_1, \,$ two subalgebras $\, A \,$ and $\, B, \,$ generated respectively by $\, a \,$ and $\, b, \,$ are isomorphic if and only if $\, \max (a) = \max (b) = l \,$, \vspace{.5mm} that is if and only if, whatever \vspace{.5mm} are $\, 0 < m < p, \,$ with $\, p \,$ prime number and $\, \dfrac mp \leq l, \,$ we have that $\, \eta_{m,p}(a) \,$ is archimedean if and only if $\, \eta_{m,p}(b) \,$ is \vspace{1.5mm} archimedean. \\
\indent Let $\, f, \, g \in Free_1, \,$ with $\, f(0)=g(0)=0. \,$ Then there is a sequence \\
$\, a_0=0<a_1<...<a_{k-1}<a_k=1 \,$ of rational numbers in $\, [0,1], \,$ together with linear functions with integer coefficients $\, p_1, \dots , p_k, \, q_1, \dots , q_k : \mathbb{R} \rightarrow \mathbb{R} \,$ such \vspace{1.5mm} that: \\
\indent 1)\ over each interval $\, [a_{i-1},a_i], \ \, f \,$ coincides with $\, p_i, \,$ and $\, g \,$ coincides \\
\hspace*{9mm} with \vspace{1mm} $\, q_i, \ \, (i = 1, \dots , k) \,$ \\
\indent 2)\ for each $\, j = 2, \dots , k, \,$  either $\, p_{j-1} \,$ is distinct from $\, p_j \,$  or else $\, q_{j-1} \,$ \\
\hspace*{9mm} is distinct from \vspace{1mm} $\, q_j.$ \\
Let the function $\, c = (f, g) :  [0, \, 1] \rightarrow [0, \, 1]^2 \,$ defined by $\, c(t) = (f(t),g(t)); \,$ the shape of range of $\, (f, g) \,$ is of course the broken line in $\, [0, \, 1]^2 \,$ joining, in the order, the \vspace{.7mm} points \\
\hspace*{10mm} $P_0 \equiv (f(a_0),g(a_0)), \ P_1 \equiv (f(a_1),g(a_1)), \dots , \vspace{.5mm} P_k \equiv (f(a_k),g(a_k)). \,$ \\
We will name these points, that are known as the nodes of the range of $\, (f, g), \,$ \textit{extremals} of the pair $\, f \,$ and \vspace{1.5mm} $\, g.$  \\
\textsc{Proposition 1.4.} - \it If $\, f, \, g \,$ and $\, f_1, \, g_1 \,$ are two pairs of $\, Free_1 \,$ elements such that the range of $\, (f, g) \,$ coincides with the range of $\, (f_1, g_1), \,$ then the algebra generated by $\, f, \, g \,$ is isomorphic to the algebra generated by \vspace{1.5mm} $\, f_1, \, g_1.$ \rm \\
\textsc{Proof.} - It follows straight from \vspace{2mm} proposition 1.3. \hspace{\stretch{1}} $\blacktriangle$ \\
\textsc{Example 1.5.} - The functions
\small
$$f(x)=\left\{
\begin{array}{lcl}
6x & \textrm{for} & 0\leq x\leq \dfrac 16 \vspace{1.5mm} \\
3-12x & \textrm{for} & \dfrac 16 \leq x\leq \dfrac 14 \vspace{1.5mm} \\
12x-3 & \textrm{for} & \dfrac 14 \leq x\leq \dfrac 13 \vspace{1.5mm} \\
1 & \textrm{for} & \dfrac 13\leq x\leq 1
\end{array}
\right. \textrm{and} \ \ g(x)=\left\{
\begin{array}{lcl}
4x & \textrm{for} & 0\leq x\leq \dfrac 16 \vspace{1.5mm} \\
2-8x & \textrm{for} & \dfrac 16 \leq x\leq \dfrac 14 \vspace{1.5mm} \\
8x-2 & \textrm{for} & \dfrac 14 \leq x\leq \dfrac 38 \vspace{1.5mm} \\
1 & \textrm{for} & \dfrac 38 \leq x\leq 1 \vspace{1.5mm}
\end{array}
\right. $$
\normalsize
and the \vspace{1.5mm} functions
$$f_1(x)=\left\{
\begin{array}{lcl}
3x & \textrm{for} & 0\leq x\leq \dfrac 13 \\
1 & \textrm{for} & \dfrac 13\leq x\leq 1
\end{array}
\right. \ \textrm{and} \quad g_1(x)=\left\{
\begin{array}{lcl}
2x & \textrm{for} & 0\leq x\leq \dfrac 12 \\
1 & \textrm{for} & \dfrac 12\leq x\leq 1 \vspace{1.5mm}
\end{array}
\right. $$
generate isomorphic \vspace{1mm} algebras.
\begin{center}
\includegraphics[scale =.3]{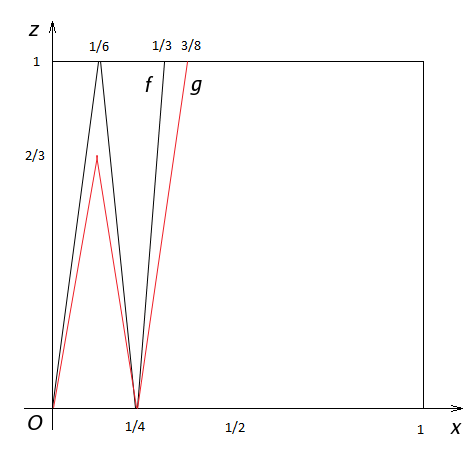}
\hspace{2mm}
\includegraphics[scale =.3]{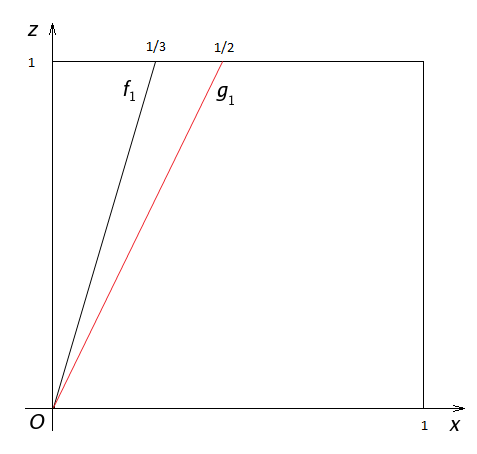} \\
Fig. \vspace{2mm} 4
\end{center}
\textsc{Corollary 1.6.} - \it If  $\, f, \, g \,$ and $\, f_1, \, g_1 \,$ are two pairs of $\, Free_1 \,$ elements \\
with the same extremals, then the algebra generated by $\, f, \, g \,$ is isomorphic \\
to the algebra generated by \vspace{4mm} $\, f_1, \, g_1.$ \rm \\
\textbf{2 Projective \vspace{2mm} algebras} \\
We recall that: \\
- A MV-algebra $\, A \,$ is named \emph{projective} if, given any MV-algebras $\, B \,$ and $\, C, \,$ for each epimorphism $\, g \,$ from $\, C \,$ on $\, B \,$ and each omomorphism $\, f \,$ from $\, A \,$ in $\, B, \,$ then an omomorphism $\, h \,$ from $\, A \,$ in $\, C, \,$ such that $\, h \circ g = f, \,$ always exists. \\
- The MV-algebra $\, B \,$ is said a \textit{retract} of the MV-algebra $\, A, \,$ if there are \\
a monomorphism $\, \chi : B \rightarrow A \,$ and an epimorphism $\, \varepsilon : A \rightarrow B \,$ such that $\, \varepsilon \chi : B \rightarrow B \,$ is the identity; \qquad so if $\, B \,$ is a retract of $\, A \,$, then $\, \chi (B) \,$ is \\
a subalgebra of $\, A \,$ such as the endomorphism $\, \chi \varepsilon \,$ of $\, A \,$ restricted to $\, \chi(B) \,$ is the identity. \quad It is known that a retract of a free MV-algebra is \vspace{1.5mm} projective. \\
\indent We try, in this paragraph, to characterize the projective subalgebras \\
of $\, Free_n. \,$ \quad First of all, we observe, that if $\, A \,$ is a projective MV-algebra \\
essentially $\, n$-generated then it is a retract of $\, Free_n. \,$ \quad More generally, \\
if $\, A, \ n$-generated, is a retract of $\, Free_m \,$ with $\, m>n, \,$ then $\, A \,$ is a retract \\
of \vspace{1.5mm} $\, Free_n$. \\
\textsc{Lemma 2.1.} - \it If $\, x_1 , \dots , x_n , \ \, u_1 , \dots , u_n \,$ are elements of a MV-algebra, $\, t \,$ is \\
a term and $\, \delta(x,u) \,$ is the Chang's distance, then there is $\, k \in \mathbb{N} \,$ such that \\
\hspace*{9mm} $\, \delta\big(t(x_1 , \dots , x_n),t(u_1 , \dots , u_n)\big) \leq \vspace{1.5mm} k\big[\delta(x_1,u_1)+ \dots + \delta(x_n,u_n)\big].$ \rm \\
\textsc{Proof.} - By an easy induction from a theorem on elementary properties \\
of the Chang's \vspace{2mm} distance [2, 3.14]. \hspace{\stretch{1}} $\blacktriangle$ \\
\textsc{Proposition 2.2.} - \it Let $\, B \,$ be a $n$-generated subalgebra of the free algebra $\, A = \, Free_n \,$ \ and let be $\, g_1 \, , \dots, \, g_n \,$ the generators of $\, A. \,$ If $\, \varphi \,$ is an epimorphism from $\, A \,$ on $\, B, \,$ then $\, \varphi(g_1)=d_1 \, , \dots , \, \varphi(g_n)=d_n \,$ generate $\, B \,$ and the kernel of $\, \varphi \,$ is the ideal $\, I \,$ generated by \vspace{1mm} $\, \delta(d_1,g_1) \, , \dots , \, \delta (d_n,g_n).$ \rm \\
\textsc{Proof.} - If $\, a \in A, \,$ where $\, a=t(g_1, \dots , g_n), \,$ is such that $\, \varphi(a)=0, \,$ then: \\
$a = \delta(a,0) = \delta\big(t(g_1, \dots , g_n),t(d_1, \dots , d_n)\big) \leq k\big[\delta(g_1,d_1) + \dots + \delta(g_n,d_n)\big],$ for an appropriate $\, k\in \mathbb{N}, \,$ that \vspace{2mm} is $\, a\in I.$ \hspace{\stretch{1}} $\blacktriangle$ \\
Thus, in order that $\, \varphi \,$ restricted to $\, B \,$ is an injection, it is necessary (and sufficient) that if $\, b\in B \,$ and $\, b \leq k\big[\delta(g_1,d_1)+ \dots + \delta(g_n,d_n) \big], \,$ then \vspace{2mm} $\, b=0.$ \\
Finally with the notations and conditions set above (proposition 2.2), we introduce a set $\, K, \,$ that we say ``\textit{equalizer}'' of \vspace{2mm} $\, B:$ \\
\textsc{Definition 2.3.} - \it $\, K=\big\{\textbf{x}\in [0,1]^n \, \pmb{:} \, d_1(\textbf{x}) = g_1(\textbf{x}) \, , \dots , \, d_n(\textbf{x}) = g_n(\textbf{x}) \big\}.$ \vspace{3mm} \\
\textsc{Theorem 2.4.} - \it Let $\, B \,$ be a $n$-generated subalgebra of $\, A = Free_n, \ \ \varphi$ \\
an epimorphism from $\, A \,$ on $\, B \,$ and let $\, g_1 \, , \dots, \, g_n \,$ be the generators of $\, A. \,$ \\
Set $\, d_1 = \varphi(g_1) \, , \dots , \, d_n = \varphi(g_n), \,$ then the epimorphism $\, \varphi \,$ restricted to $\, B, \,$ is an isomorphism if and only if for every $\, \textbf{u} \in [0,1]^n \,$ there is $\, \textbf{x} \in K \,$ \\
such that \\
\hspace*{35mm} $\, d_1(\textbf{u}) = d_1(\textbf{x}) \, , \dots , \, d_n(\textbf{u}) = \vspace{3mm} d_n(\textbf{x}).$ \rm \\
\textsc{Proof.} - If for $\, b\in B, \ b=t(d_1, \dots ,d_n), \,$ were \\
$\, b\leq k\big[\delta (d_1,g_1)+ \dots + \delta (d_n,g_n)\big] \,$ and $\, b\neq 0, \,$ then there would be \vspace{1mm} $\, \textbf{\textit{u}} \in [0,1]^n $ \\
such that: $\, b(\textbf{\textit{u}}) = t\big(d_1(\textbf{\textit{u}}), \dots ,d_n(\textbf{\textit{u}})\big)\neq 0; \quad $ but by \vspace{1mm} hypothesis, \\
$\, \textbf{\textit{x}}_0 \in K $ exists such \vspace{1mm} that $\, d_1(\textbf{\textit{u}}) = d_1(\textbf{\textit{x}}_0) \, , \dots , \, d_n(\textbf{\textit{u}}) = d_n(\textbf{\textit{x}}_0) \,$ \\
and so \vspace{1mm} $\, t\big(d_1(\textbf{\textit{u}}), \dots ,d_n(\textbf{\textit{u}})\big) = t\big(d_1(\textbf{\textit{x}}_0), \dots ,d_n(\textbf{\textit{x}}_0)\big) = b(\textbf{\textit{x}}_0)\neq 0.$ \\
But, since $\, d_1(\textbf{\textit{x}}_0) = g_1(\textbf{\textit{x}}_0) \, , \dots , \, d_n(\textbf{\textit{x}}_0) = g_n(\textbf{\textit{x}}_0) \,$ we \vspace{1mm} have: \\
$\big(k\big[\delta (d_1,g_1)+ \dots + \delta (d_n,g_n)\big]\big)(\textbf{\textit{x}}_0) \vspace{1mm} = \\
\hspace*{28mm} k\big[\delta \big(d_1(\textbf{\textit{x}}_0),g_1(\textbf{\textit{x}}_0)\big) + \dots + \delta\big(d_n(\textbf{\textit{x}}_0),g_n(\textbf{\textit{x}}_0) \big) \big] = \vspace{1mm} 0$ \\
and so would be \vspace{1mm} $\, 0\neq b(\textbf{\textit{x}}_0)\leq \big(k\big[\delta (d_1,g_1)+ \dots + \delta (d_n,g_n)\big]\big)(\textbf{\textit{x}}_0)=0.$ \\
Vice versa, let $\, \textbf{\textit{u}}_1 \in [0,1]^n \,$ be such that for every $\, \textbf{\textit{x}} \in K,$ \\
\hspace*{25mm} $d_i(\textbf{\textit{u}}_1) \neq d_i(\textbf{\textit{x}}) \qquad \textrm{for some} \ \ \vspace{.5mm} i, \ 1 \leq i \leq n, $ \\
consequently there is $\, \varepsilon > 0 \,$ such as for every $\, \textbf{\textit{v}} \in \mathscr{I}(\textbf{\textit{u}}_1, \varepsilon ) \,$ and $\, \textbf{\textit{x}} \in K \,$ \\
it is \ $d_i(\textbf{\textit{v}}) \neq d_i(\textbf{\textit{x}}). \quad$ In particular there are \\
integers $\, 0 < m_1 < p_1, \, \dots \, ,0 < m_n < p_n,$ with $\, p_1, \dots ,p_2 \,$ prime numbers, $\, \textbf{\textit{u}}_2 \in \mathscr{I} (\textbf{\textit{u}}_1, \varepsilon ) \,$ and $\, \varrho > 0 \,$ such that:
$$d_1(\textbf{\textit{u}}_2) = \dfrac{m_1}{p_1}, \, \dots \, ,d_n(\textbf{\textit{u}}_2) = \dfrac{m_n}{p_n} \ \ \textrm{and} \ \ \mathscr{I} (\textbf{\textit{u}}_2, \varrho ) \subseteq \mathscr{I} (\textbf{\textit{u}}_1, \varepsilon ).$$
Then we can choose $\, h \in \mathbb{N} \,$ so that \vspace{2mm} the element \\
\hspace*{45mm} $\displaystyle b = \Big(h \bigvee_{i=1}^n \eta_{m_i,p_i}(d_i) \vspace{2mm} \Big)'$ \\
is nonzero only in $\, \mathscr{I} (\textbf{\textit{u}}_2, \varrho ). \,$ Thus, for $\, k \in \mathbb{N} \,$ large enough, \\
is $\, b \leq k\big[\delta(g_1,d_1)+ \dots + \delta(g_n,d_n) \big], \,$ but $\, b \neq 0 \,$ \vspace{3mm} also. \hspace{\stretch{1}} $\blacktriangle$ \\
\textsc{Corollary 2.5.} - \it Let $\, A \,$ be a $n$-generated MV-algebra. $\, A \,$ is projective \\
if and only if is isomorphic to a subalgebra $\, B \,$ of $\, Free_n, \,$ whereby, \\
however we choose $\, \textbf{u} \in [0,1]^n \,$ there is $\, \textbf{x} \in K \,$ such that \\
\hspace*{35mm} $\, d_1(\textbf{u}) = d_1(\textbf{x}) \, , \dots , \, d_n(\textbf{u}) = \vspace{3mm} d_n(\textbf{x}).$ \\ \rm
\textsc{Corollary 2.6 (Di Nola).} - \it Every monogenerated subalgebra of $\, Free_1 \,$ \\
is \vspace{1.5mm} projective. \rm

We devote the next section to some constructive applications \\
of theorem 2.4 for bigenerated projective \vspace{4mm} algebras. \\
\textbf{3  Bigenerated \vspace{2mm} algebras} \\
Since the equalizer $\, K \,$ must be a finite union of triangular simplexes \\
$\, T_1, \, T_2, \, ..., \,T_k \,$ and each triangular simplex $\, T_i \,$ is determinated by a system
$$ \left\{
\begin{array}{l}
a_{i\phantom{1}} x + b_{i\phantom{1}} y + c_{i\phantom{1}} \geq 0 \\
a_{i1} x + b_{i1} y + c_{i1} \geq 0  \hspace{32mm} (\alpha_i) \\
a_{i2} x + b_{i2} y + c_{i2} \geq 0
\end{array}
\right. $$
of three linear inequalities in two variables, so we can state, according \\
to previous \vspace{1.5mm} observations: \\
\textsc{Theorem 3.1.} - \it A subalgebra $\, B \,$ of $\, Free_2 \,$ generated by $\, a \,$ and \vspace{1.5mm} $\, b,$ \\
is projective if and only if the pairs $\, \big(\dfrac{m_1}{p_1}, \dfrac{m_2}{p_2}\big), \,$ with \vspace{1.2mm} $\, 0<m_j<p_j \,$ \\
and $\, p_j \,$ prime $\, (j=1,2), \,$ that satisfy one of the systems $\, (\alpha _1) , (\alpha _2) , \dots , (\alpha _k), \,$ are \underline{all} and \underline{onl}y those for wich the element $\, \eta_{m_1,p_1}(a) \vee \eta_{m_2,p_2}(b) \,$ is not archi- medean. \vspace{1.5mm} \rm \\
\textsc{Proof.} - By \vspace{1.5mm} corollary 2.5. \hspace{\stretch{1}} $\blacktriangle$ \\
\indent Therefore, it is important to establish wich are the triangular simplexes that give rise the projective algebras. Certainly, the union of such simplexes must be connected and it must contain $\, \textbf{0} \,$. But it can be not convex, as it is easy to see. \\
\indent Now, we point out a geometrical fact quite elementary, but basic for the \vspace{1.5mm} following. \\
\textsc{Remark 3.2.} - On the plane $\, Oxz, \,$ we consider the points \\
\hspace*{20mm} $O \equiv (0,0), \ P \equiv (a,h), \ K \equiv (0,k), \ Q \equiv (a,l),$ \\
\hspace*{20mm} $O' \equiv (b,0), \ P' \equiv (c,h), \ K' \equiv (b,k) \ \textrm{and} \ Q' \equiv (c,l),$ \\
where the real numbers $\, a, \; b \ \textrm{and} \ c, \,$ each nonzero, are pairwise different. Similarly, $\, h, \; k \ \textrm{and} \ l \,$ are supposed to be all nonzero and pairwise different. \ If $\, s \; \textrm{and} \; t \in \mathbb{R}, \,$ the line $\, x = s \,$ intersects the line through $\, O \, P \,$ and the line through $\, K \, Q \,$ respectively at points $\, S \equiv (s,u) \ \textrm{and} \ S' \equiv (s,u'); \,$ and the line $\, x = t \,$ intersects the line through $\, O' \, P' \,$ and the line through $\, K' \, Q' \,$ respectively at points $\, T \equiv (t,v) \ \textrm{and} \ T' \equiv (t,v'). \ \ $ \it \\
If, for a given $\, s, \,$ we choose $\, t = b+[s(c-b)/a] \,$ $($or if, for a given $\, t, \,$ \\
we choose $\, s = a(t-b)/(c-b)), \,$ then $\, u = v \ \textrm{and} \ u'=v'.$ \rm
\begin{center}
\includegraphics[scale =.6]{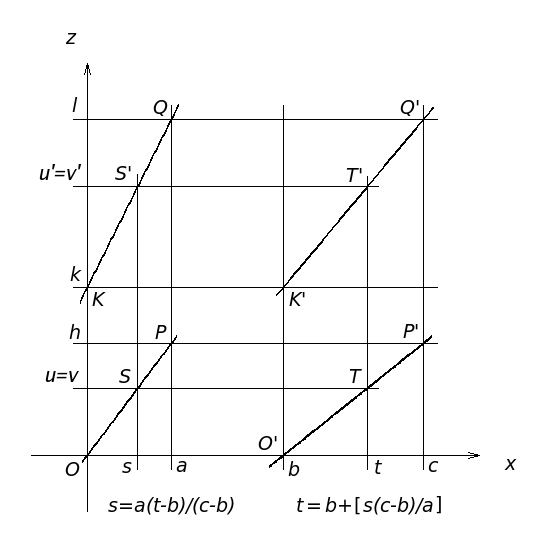} \vspace{-2mm} \\
Fig. 5
\end{center}
\newpage
In the space $\, \mathbb{R}^3, \,$ we denote by $\, \Delta_{\mathcal{U}} \,$ the triangle 
whose vertices belong to the set $\qquad \mathcal{U} = \big\{(t_1,\, v_1,\, z_1),\; (t_2,\, 
v_2,\, z_2),\; (t_3,\, v_3,\, z_3) \big\},$ \qquad where \\
the points $\, (t_1,\, v_1,\, 0),\; (t_2,\, v_2,\, 0),\; (t_3,\, v_3,\, 0) 
\,$ are not aligned. We denote by $\, z = f_{\mathcal{U}}(t,\, v) \,$ the linear 
function corresponding to the plane settled by the points of \vspace{2mm} $\, \mathcal{U}.$ \\
\textsc{Proposition 3.3.} \it Let the \vspace{1mm} sets: \\
\begin{tabular}{ll}
\hspace{25mm} & $\mathcal{A} = \big\{(x_1,\, y_1,\, a_1),\; (x_2,\, y_2,\, a_2),\; (x_3,\, y_3,\, 
a_3) \big\},$ \\
& $\mathcal{B} = \big\{(x_1,\, y_1,\, b_1),\; (x_2,\, y_2,\, b_2),\; (x_3,\, y_3,\, b_3) \big\},$ \\
& $\mathcal{C} = \big\{(x_1,\, y_1,\, 0),\; (x_2,\, y_2,\, 0),\; (x_3,\, y_3,\, 0) \big\},$ 
\vspace{1mm}
\end{tabular} \\
with the condition that the points of $\, \mathcal{C} \,$ are not aligned. \\
Arbitrarily chosen \\
a set $\mathcal{J} = \big\{(\xi_1,\, \eta_1,\, 0),\; (\xi_2,\, \eta_2,\, 0),\; (\xi_3,\, 
\eta_3,\, 0) \big\},$ whose points are not aligned, a 3-list 
$\, (i_1,\, i_2,\, i_3) \,$ in $\, \{1,\, 2,\, 3\},$ \\
and consequently also the \vspace{1mm} sets: \\
\begin{tabular}{ll}
& $\,\ \mathcal{L} = \big\{(\xi_1,\, \eta_1,\, a_{i_1}),\; (\xi_2,\, \eta_2,\, a_{i_2}),\; 
(\xi_3,\, \eta_3,\, a_{i_3}) \big\}$ and \\
\hspace{19mm} & $\mathcal{M} = \big\{(\xi_1,\, \eta_1,\, b_{i_1}),\; (\xi_2,\, \eta_2,\, 
b_{i_2}),\; (\xi_3,\, \eta_3,\, b_{i_3}) \big\}.$ \vspace{1mm}
\end{tabular} \\
Then, for every $\, (\xi_0,\, \eta_0,\, 0) \in \Delta_{\mathcal{J}} \,$ there is $\, (x_0,\, 
y_0,\, 0) \in \Delta_{\mathcal{C}} \,$ such that \\
\hspace*{22mm} $\, f_{\mathcal{L}}(\xi_0,\, \eta_0) = f_{\mathcal{A}}(x_0,\, y_0) \ \ $ and $\ \
f_{\mathcal{M}}(\xi_0,\, \eta_0) = f_{\mathcal{B}}(x_0,\, y_0).$ \vspace{1mm} \rm \\
\textsc{Proof.} - Suppose first, that  $\, (i_1,\, i_2,\, i_3) \,$ is a permutation of $\, \{1,\, 
2,\, 3\}. \,$ \\
Let $\, \alpha \,$ the plane passing through the points $\, (\xi_1,\, \eta_1,\, a_{i_1}),\; 
(\xi_1,\, \eta_1,\, b_{i_1}) \,$ \\
and $\, (\xi_0,\, \eta_0,\, 0). \,$ This plan meets the plane $\, \beta, \,$ passing through the points \\
$\, (\xi_2,\, \eta_2,\, a_{i_2}),\; (\xi_2,\, \eta_2,\, b_{i_2}) \,$ and $\, (\xi_3,\, \eta_3,\, a_{i_3}), \,$ along a straight line $\, r \,$ \\
which intersects the side of the triangle $\, \Delta_{\mathcal{L}}, \,$ of extremes 
$\, (\xi_2,\, \eta_2,\, a_{i_2}) \,$ and $\, (\xi_3,\, \eta_3,\, a_{i_3}) \,$ and the side of the 
triangle $\, \Delta_{\mathcal{M}}, \,$ of extremes $\, (\xi_2,\, \eta_2,\, b_{i_2}) \,$ and \\
$\, (\xi_3,\, \eta_3,\, b_{i_3}), \,$ respectively, in points $\, (\xi_4,\, \eta_4,\, h) \,$ and 
$\, (\xi_4,\, \eta_4,\, k) \,$ \hspace{\stretch{1}} (Fig. 6.0).\\
For the remark 3.2, on the side of the triangle $\, \Delta_{\mathcal{C}}, \,$ of extremes $\, 
(x_{i_2},\, y_{i_2},\, 0) \,$ \\
and $\, (x_{i_3},\, y_{i_3},\, 0), \,$ there is a point $\, (x_4,\, y_4,\, 0) \,$ such that $\, 
f_{\mathcal{A}}(x_4,\, y_4) = h \,$ \\
and $\, f_{\mathcal{B}}(x_4,\, y_4) = k \,$ (Fig 6.1). Just now consider the plane $\, \gamma \,$ 
passing through the points $\, (x_{i_1},\, y_{i_1},\, a_{i_1}),\; (x_{i_1},\, y_{i_1},\, b_{i_1}) \,$ e $\, 
(x_4,\, y_4,\, 0) \,$ and the assertion follows by a further application of the remark 3.2 \quad (Fig. 6.2). \\
If $\, (i_1,\, i_2,\, i_3) \,$ is not a permutation of  $\, \{1,\, 2,\, 3\}, \,$ similar 
proofs \vspace{2mm} hold. \hspace{\stretch{1}} $\blacktriangle$
\begin{center}
\includegraphics[scale =.8]{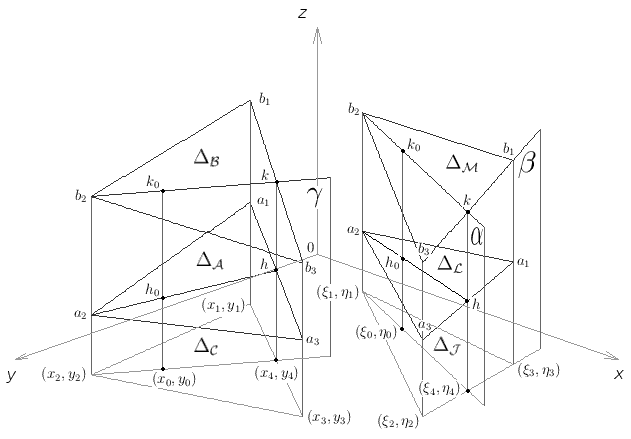} \vspace{1mm} \\
\footnotesize Here the 3-list is (2,3,1). \quad\ For simplicity, points in plane $z = 0$ are \vspace{-1mm} denoted \\
with the first two coordinates, while the other points only with their \vspace{1mm} heights. \normalsize \\
Fig. \vspace{6mm} 6.0 \\
\includegraphics[scale =.8]{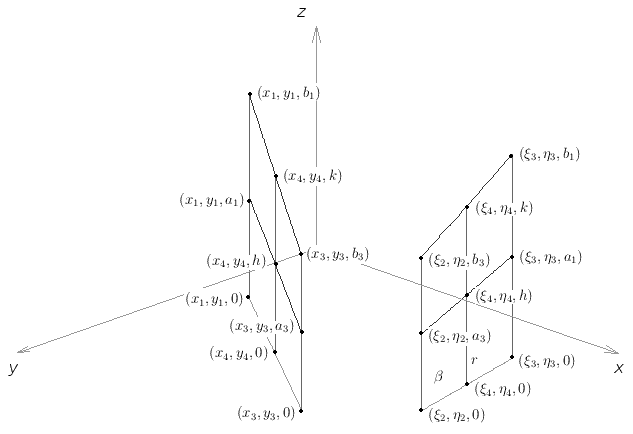} \vspace{1mm} \\
\footnotesize First application of remark 3.2, this \vspace{-1mm} yields the point \\
$\, (x_4,\, y_4,\, 0) \,$ \ \ corresponding to the point \ \ \vspace{1mm} $\, (\xi_4,\, \eta_4,\,0) \,$ \normalsize \\
Fig. 6.1
\newpage \noindent
\includegraphics[scale =.8]{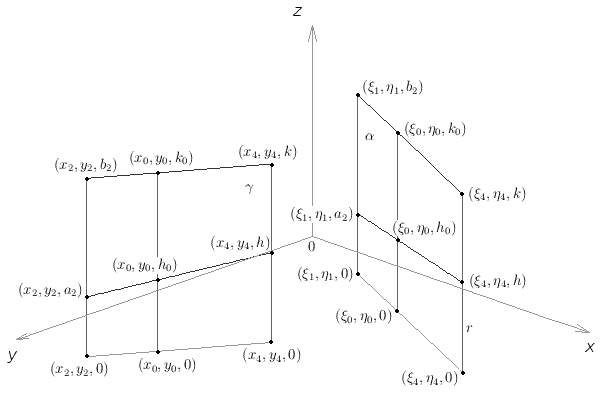} \vspace{1mm} \\
\footnotesize Second application of remark 3.2, this \vspace{-1mm} yields the point \\
$\, (x_0,\, y_0,\, 0) \,$ \quad corresponding to the point \quad \vspace{1mm} $\, (\xi_0,\, \eta_0,\, 0) \,$ \normalsize \\
Fig. \vspace{3mm} 6.2
\end{center}
We conclude finally exposing some \vspace{2mm} applications. \\
$i)$ - A first very simple case arises if the equalizer is:
$$K = \big\{(x,y) \in [0,1]^2 \, : \, g_1(y) \leq x \leq g_2(y) \ \; \textbf{\textrm{and}} \ \; f_1(x) \leq y \leq f_2(x) \big\},$$
where $\, f_1(x), \ f_2(x), \ g_1(y) \ \textrm{and} \ g_2(y) \,$ are McNaughton functions \\
and the follow conditions are satisfied:
$$\, f_1(x) \leq f_2(x), \ \; g_1(y) \leq \vspace{-6.5mm} g_2(y); $$
\hspace*{\stretch{1}} \vspace{-4mm} $(\bullet)$
$$y \geq f_1 \big(g_1(y) \big), \ \; y \leq f_2 \big(g_1(y) \big), \ \; y \leq f_2 \big(g_2(y) \big), \ \; y \vspace{2mm} \geq f_1 \big(g_2(y) \big).$$
Set\,: \\
\begin{minipage}[c]{70mm}
\begin{center}
$A = \big\{(x,y) \in [0,1]^2 \, : \, y \leq f_1(x) \vspace{1.5mm} \big\},$ \\
$B = \big\{(x,y) \in [0,1]^2 \, : \, y \geq f_2(x) \vspace{1.5mm} \big\},$ \\
$C = \big\{(x,y) \in [0,1]^2 \, : \, x \leq g_1(y) \vspace{1.5mm} \big\},$ \\
$D = \big\{(x,y) \in [0,1]^2 \, : \, x \geq g_2(y) \big\}.$
\end{center}
\end{minipage} \hspace{3mm}
\begin{minipage}[c]{60mm}
\vspace*{4mm}
\hspace*{5mm}\includegraphics[scale =.3]{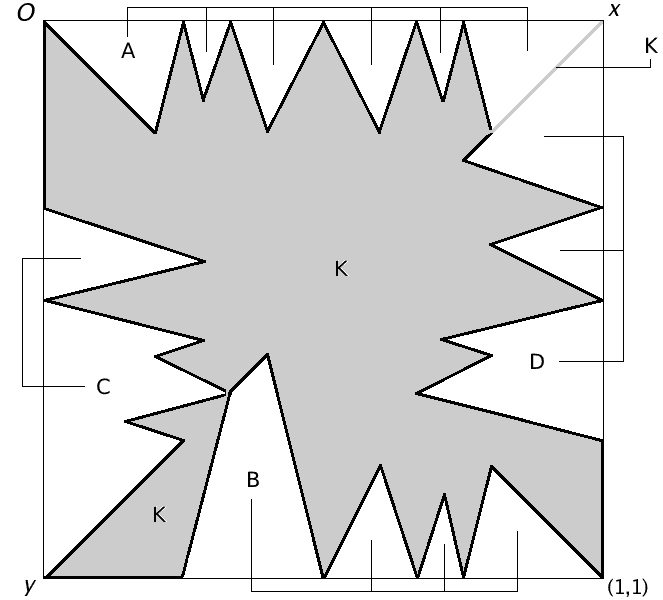} \vspace{2mm} \\
\hspace*{22mm} \vspace{4mm} Fig. 7
\end{minipage} \\
\normalsize
The functions
$$d_1(x,y)=\left\{
\begin{array}{lll}
x & \textrm{if} & (x,y) \in A \cup K \cup B, \\
g_1(y) \ \; & \, \raisebox{-1.5mm}{''} & (x,y) \in C, \\
g_2(y) & \, \raisebox{-1.5mm}{''} & (x,y) \in D
\end{array}
\right. $$
$$d_2(x,y)=\left\{
\begin{array}{lll}
y & \textrm{if} & (x,y) \in C \cup K \cup D, \\
g_1(y) \ \; & \, \raisebox{-1.5mm}{''} & (x,y) \in A, \\
g_2(y) & \, \raisebox{-1.5mm}{''} & (x,y) \in B
\end{array}
\right. \vspace{3mm} $$
generate, by the corollary 2.5 and the proposition 3.3, a projective \vspace{.5mm} subalgebra of \vspace{1.5mm} $\, Free_2.$ \\
\textsc{Example 3.4.} - The McNaughton functions
\small
$$f_1(x)=\left\{
\begin{array}{lcl}
0 & \textrm{for} & 0\leq x\leq \dfrac 13 \vspace{1.5mm} \\
1-2x & \textrm{for} & \dfrac 13 \leq x\leq \dfrac 12 \vspace{1.5mm} \\
2x-1 & \textrm{for} & \dfrac 12 \leq x\leq \dfrac 23 \vspace{1.5mm} \\
1-x & \textrm{for} & \dfrac 23\leq x\leq 1
\end{array}
\right. \ , \quad f_2(x)=\left\{
\begin{array}{lcl}
1-x & \textrm{for} & 0 \leq x\leq \dfrac 12 \vspace{1.5mm} \\
x & \textrm{for} & \dfrac 12 \leq x\leq 1
\end{array}
\right. \ ;$$
$$g_1(y)=\left\{
\begin{array}{lcl}
0 & \textrm{for} & 0 \leq x\leq \dfrac 14 \vspace{1.5mm} \\
4y-1 & \textrm{for} & \dfrac 14 \leq y\leq \dfrac 13 \vspace{1.5mm} \\
1-2y & \textrm{for} & \dfrac 13 \leq y\leq \dfrac 12 \vspace{1.5mm} \\
2y-1 & \textrm{for} & \dfrac 12 \leq y\leq \dfrac 23 \vspace{1.5mm} \\
1-y & \textrm{for} & \dfrac 23\leq y\leq 1
\end{array}
\right. \ , \quad g_2(y)=\left\{
\begin{array}{lcl}
1 & \textrm{for} & 0 \leq y\leq \dfrac 13 \vspace{1.5mm} \\
2-3y & \textrm{for} & \dfrac 13 \leq y\leq \dfrac 25 \vspace{1.5mm} \\
2y & \textrm{for} & \dfrac 25 \leq y\leq \dfrac 12 \vspace{1.5mm} \\
1 & \textrm{for} & \dfrac 12 \leq y\leq \dfrac 23 \vspace{1.5mm} \\
3-3y & \textrm{for} & \dfrac 23 \leq y\leq \dfrac 57 \vspace{1.5mm} \\
4y-2 & \textrm{for} & \dfrac 57 \leq y\leq \dfrac 34 \vspace{1.5mm} \\
1 & \textrm{for} & \dfrac 34\leq y\leq 1
\end{array}
\right. \ \,$$
\normalsize
satisfy the conditions $\, (\bullet);$
\begin{center}
\includegraphics[scale =.30]{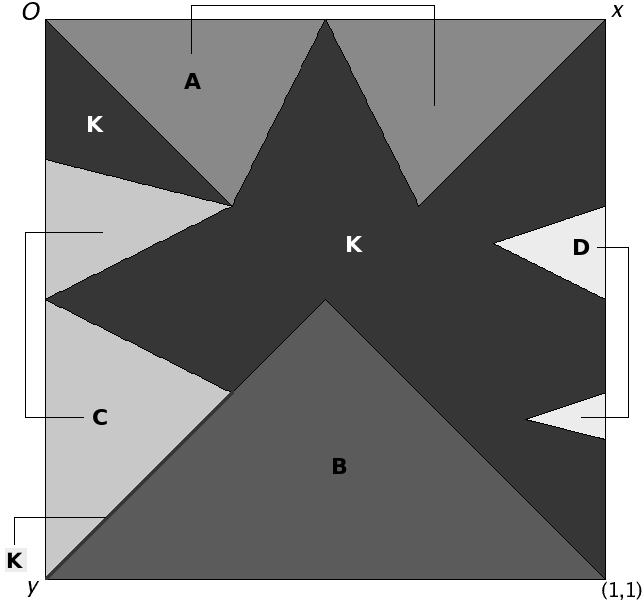} \\
Fig. 8
\end{center}
Therefore
\small
$$d_1(x,y)=\left\{
\begin{array}{lll}
4y-1 \ & \textrm{if} & x \leq 4y-1 \ \; \textrm{\textbf{and}} \ \; y \leq 1/3, \\
1-2y & \, \raisebox{-1.5mm}{''} & x \leq 1-2y \ \; \textrm{\textbf{and}} \ \; y \geq 1/3, \\
2y-1 & \, \raisebox{-1.5mm}{''} & x \leq 2y-1 \ \; \textrm{\textbf{and}} \ \; y \leq 2/3, \\
1-y & \, \raisebox{-1.5mm}{''} & x \leq 1-y \ \; \textrm{\textbf{and}} \ \; y \geq 2/3, \\
2-3y & \, \raisebox{-1.5mm}{''} & x \geq 2-3y \ \; \textrm{\textbf{and}} \ \; y \leq 2/5, \\
2y & \, \raisebox{-1.5mm}{''} & x \geq 2y \ \; \textrm{\textbf{and}} \ \; y \geq 2/5, \\
3-3y & \, \raisebox{-1.5mm}{''} & x \geq 3-3y \ \; \textrm{\textbf{and}} \ \; y \leq 5/7, \\
4y-2 & \, \raisebox{-1.5mm}{''} & x \geq 4y-2 \ \; \textrm{\textbf{and}} \ \; y \geq 5/7, \\
x & & \textrm{otherwise}
\end{array}
\right. $$
\normalsize
\hspace*{25mm} and
\small
$$d_2(x,y)=\left\{
\begin{array}{lll}
x & \textrm{if} & y \leq x \ \; \textrm{\textbf{and}} \ \; x \leq 1/3, \\
1-2x \ & \, \raisebox{-1.5mm}{''} & y \leq 1-2x \ \; \textrm{\textbf{and}} \ \; x \geq 1/3, \\
2x-1 & \, \raisebox{-1.5mm}{''} & y \leq 2x-1 \ \; \textrm{\textbf{and}} \ \; x \leq 2/3, \\
1-x & \, \raisebox{-1.5mm}{''} & y \leq 1-x \ \; \textrm{\textbf{and}} \ \; x \geq 2/3, \\
1-x & \, \raisebox{-1.5mm}{''} & y \geq 1-x \ \; \textrm{\textbf{and}} \ \; x \leq 1/2, \\
x & \, \raisebox{-1.5mm}{''} & y \geq x \ \; \textrm{\textbf{and}} \ \; x \geq 1/2, \\
y & & \textrm{otherwise}
\end{array}
\right. \vspace{5mm} $$
\normalsize
are the generators of a projective subalgebra of $\, Free_2.$
\begin{center}
\includegraphics[scale =.27]{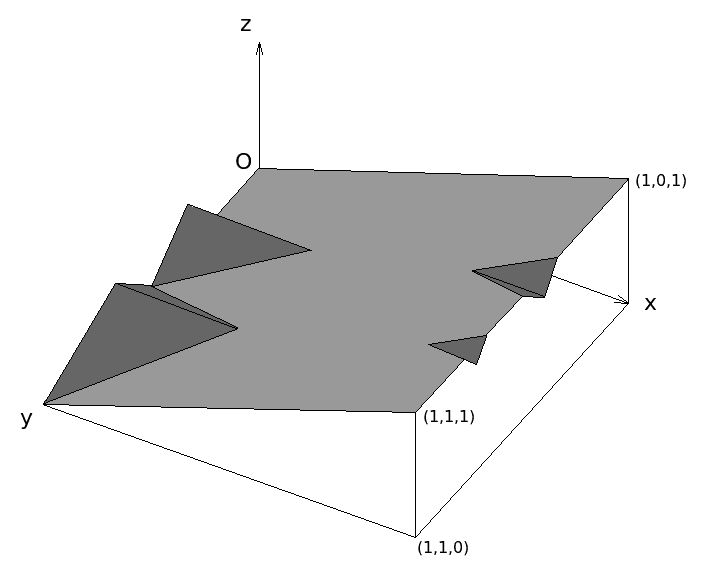}
\includegraphics[scale =.27]{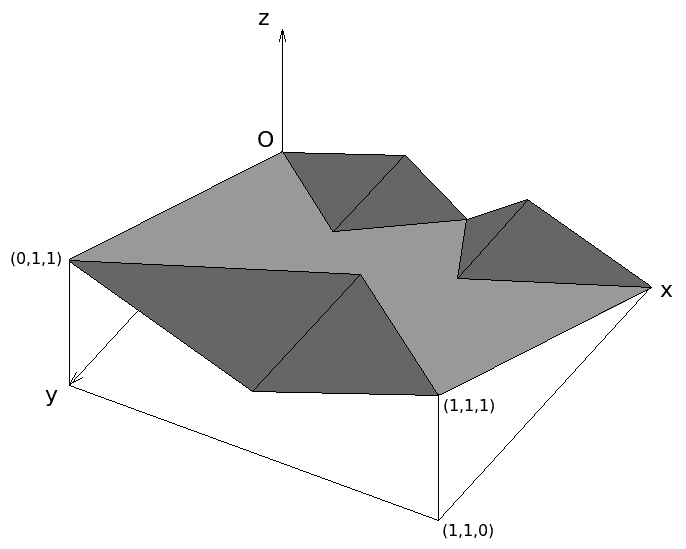} \vspace{.1mm} \\
\scriptsize
$d_1(x,y) \hspace{50mm} \vspace{2mm} d_2(x,y)$ \\
\normalsize
Fig. \vspace{2mm} 9
\end{center}
\textsc{Example 3.5.} - We consider the subalgebra $\, A \,$ of $\, Free_1 \,$ generated by the 
functions: \vspace{-2mm}
\small
$$f(x)=\left\{
\begin{array}{lcl}
3x & \textrm{for} & 0\leq x\leq \dfrac 13 \vspace{1.5mm} \\
2-3x & \textrm{for} & \dfrac 13 \leq x\leq \dfrac 23 \vspace{1.5mm} \\
3x-2 & \textrm{for} & \dfrac 23\leq x\leq 1
\end{array}
\right. \ \; \textrm{and} \quad g(x)=\left\{
\begin{array}{lcl}
3x & \textrm{for} & 0\leq x\leq \dfrac 13 \vspace{1.5mm} \\
2-3x & \textrm{for} & \dfrac 13 \leq x\leq \dfrac 12 \vspace{1.5mm} \\
3x-1 & \textrm{for} & \dfrac 12 \leq x\leq \dfrac 23 \vspace{1.5mm} \\
3x-3 & \textrm{for} & \dfrac 23 \leq x\leq 1
\end{array}
\right. .$$
\normalsize
\begin{center}
\includegraphics[scale =.26]{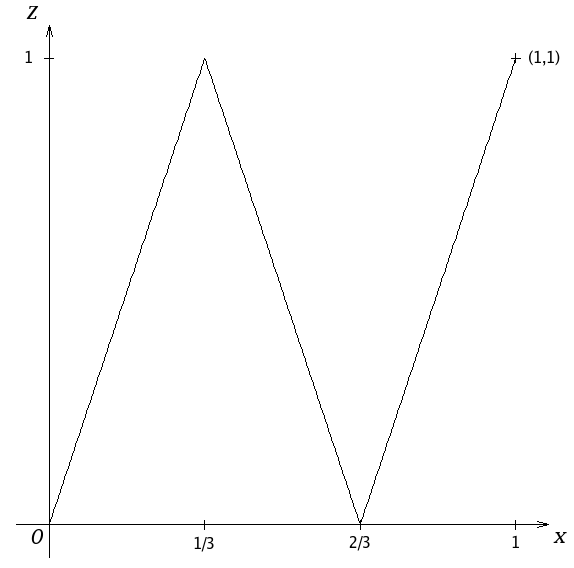} \hspace{10mm} \vspace{-1mm}
\includegraphics[scale =.26]{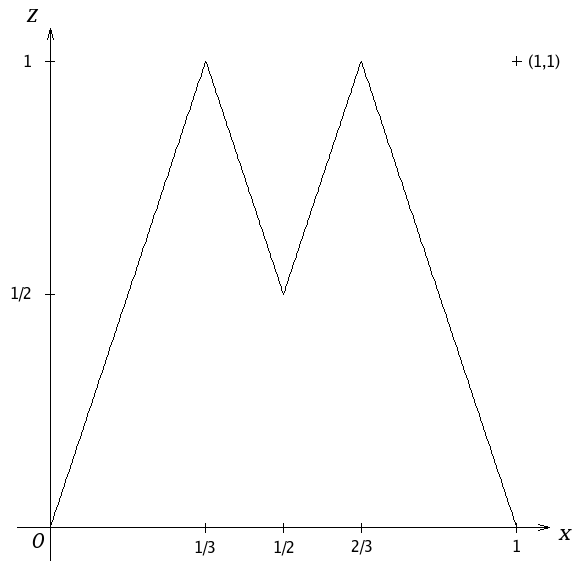} \\
\small
$f(x) \hspace{60mm} \vspace{-1mm} g(x)$ \\
\normalsize
Fig. 10
\end{center}
We can see easily that the extremals of $\, f(x) \,$ and $\, g(x) \,$ are the points
$$P_0 \equiv (0,0), \ P_1 \equiv (1,1), \ P_2 \equiv (1/2,1/2), \ P_3 \equiv (0,1)  \ \; \textrm{and} \ \: P_4 \equiv (1,0);$$
so that the range of $\, (f,g) \,$ consists in the diagonals of the square $\, [0,1]^2.$ \\
However, these diagonals may be regarded as the set of the points of $\, [0,1]^2$ \\
enclosed by the following four McNaughton functions
$$f_1(x)=\left\{
\begin{array}{lcl}
x & \textrm{for} & 0 \leq x \leq \dfrac 12 \vspace{1.5mm} \\
1-x & \textrm{for} & \dfrac 12 \leq x\leq 1
\end{array}
\right. , \quad f_2(x)=\left\{
\begin{array}{lcl}
1-x & \textrm{for} & 0 \leq x \leq \dfrac 12 \vspace{1.5mm} \\
x & \textrm{for} & \dfrac 12 \leq x\leq 1
\end{array}
\right. ,$$
$$g_1(y)=\left\{
\begin{array}{lcl}
y & \textrm{for} & 0 \leq y \leq \dfrac 12 \vspace{1.5mm} \\
1-y & \textrm{for} & \dfrac 12 \leq y \leq 1
\end{array}
\right. , \quad g_2(y)=\left\{
\begin{array}{lcl}
1-y & \textrm{for} & 0 \leq y \leq \dfrac 12 \vspace{1.5mm} \\
y & \textrm{for} & \dfrac 12 \leq y \leq 1
\end{array}
\right. ,$$
that obviously satisfy the conditions \vspace{.5mm} $\, (\bullet).$ \\
It follows that the functions
\small
$$d_1(x,y)=\left\{
\begin{array}{lll}
y \ & \textrm{if} & 0 \leq y \leq \dfrac 12 \ \; \textrm{\textbf{and}} \ \; y \geq x, \ \; \textrm{or} \\
& & \dfrac 12 \leq y \leq 1 \ \; \textrm{\textbf{and}} \ \; y \leq x, \\
1-y & \, \textrm{if} & 0 \leq y \leq \dfrac 12 \ \; \textrm{\textbf{and}} \ \; y \geq 1-x, \ \; \textrm{or} \\
& & \dfrac 12 \leq y \leq 1 \ \; \textrm{\textbf{and}} \ \; y \leq 1-x, \\
x & & \textrm{otherwise}
\end{array}
\right. $$
\normalsize
\hspace*{25mm} and
\small
$$d_2(x,y)=\left\{
\begin{array}{lll}
x \ & \textrm{if} & 0 \leq x \leq \dfrac 12 \ \; \textrm{\textbf{and}} \ \; x \geq y, \ \; \textrm{or} \\
& & \dfrac 12 \leq x \leq 1 \ \; \textrm{\textbf{and}} \ \; x \leq y, \\
1-x & \, \textrm{if} & 0 \leq x \leq \dfrac 12 \ \; \textrm{\textbf{and}} \ \; x \geq 1-y, \ \; \textrm{or} \\
& & \dfrac 12 \leq x \leq 1 \ \; \textrm{\textbf{and}} \ \; x \leq 1-y, \\
y & & \vspace{2mm} \textrm{otherwise}
\end{array}
\right. $$
\normalsize
\begin{center}
\includegraphics[scale =.32]{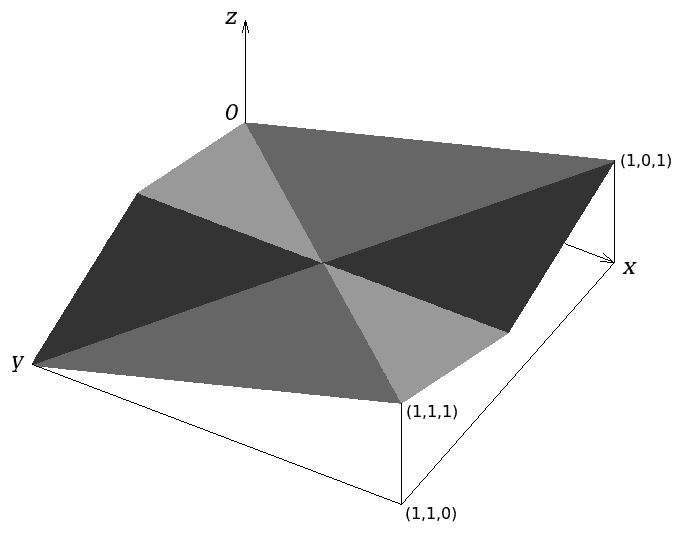}
\includegraphics[scale =.32]{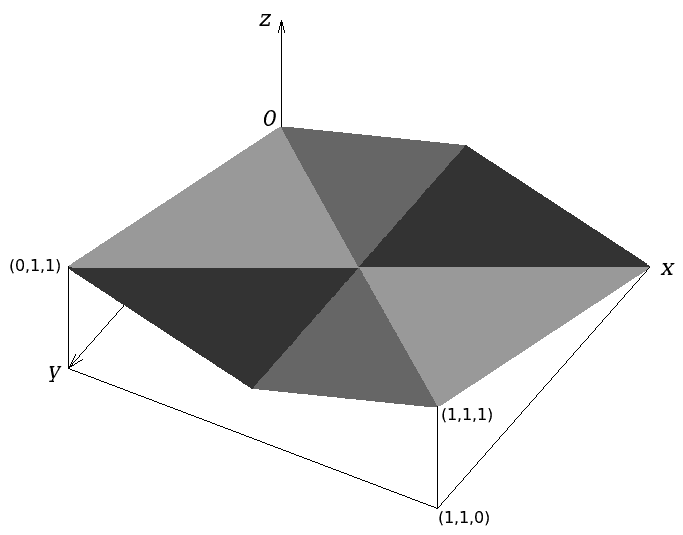} \vspace{.1mm} \\
\scriptsize
$d_1(x,y) \hspace{50mm} \vspace{2mm} d_2(x,y)$ \\
\normalsize
Fig. \vspace*{1.5mm} 11
\end{center}
are the generators of a projective subalgebra of $\, Free_2;$ \\
but this, by the proposition 1.3, is isomorphic to \vspace{2mm} $\, A.$ \\
$ii)$ - Now, let's consider the case that the equalizer is:
$$K = \big\{(x,y) \in [0,1]^2 \, : \, g_1(y) \leq x \leq g_2(y) \ \; \textbf{\textrm{and}} \ \; f_1(x) \leq y \leq f_2(x) \big\},$$
where $\, f_1(x), \ f_2(x), \ g_1(y) \ \textrm{and} \ g_2(y) \,$ are broken lines, for wich \\
the conditions:
$$\left\{
\begin{array}{lcl}
0 \leq f_1(x) \leq x & \textrm{for} & 0 \leq x \leq \dfrac 12, \vspace{1.5mm} \\
0 \leq f_1(x) \leq 1-x & \textrm{for} & \dfrac 12 \leq x \leq 1,
\end{array}
\right.$$
$$\left\{
\begin{array}{lcl}
1-x \leq f_2(x) \leq 1 & \textrm{for} & 0 \leq x \leq \dfrac 12, \vspace{1.5mm} \\
x \leq f_2(x) \leq 1 & \textrm{for} & \dfrac 12 \leq x \leq 1,
\end{array}
\right.$$
$$\left\{
\begin{array}{lcl}
0 \leq g_1(y) \leq y & \textrm{for} & 0 \leq y \leq \dfrac 12, \vspace{1.5mm} \\
0 \leq g_1(y) \leq 1-y & \textrm{for} & \dfrac 12 \leq y \leq 1,
\end{array}
\right.$$
$$\left\{
\begin{array}{lcl}
1-y \leq g_2(y) \leq 1 & \textrm{for} & 0 \leq y \leq \dfrac 12, \vspace{1.5mm} \\
y \leq g_2(y) \leq 1 & \textrm{for} & \dfrac 12 \leq y \leq 1;
\end{array}
\right. \vspace{-1.5mm} $$
are \vspace{3mm} satisfied, \\
and the sets: \\
\begin{minipage}[c]{70mm}
\begin{center}
$A =  \big\{(x,y) \in [0,1]^2 \, : \, y \leq f_1(x) \big\},$ \vspace{1.5mm} \\
$B =  \big\{(x,y) \in [0,1]^2 \, : \, y \geq f_2(x) \big\},$ \vspace{1.5mm} \\
$C =  \big\{(x,y) \in [0,1]^2 \, : \, x \leq g_1(y) \big\}.$ \vspace{1.5mm} \\
$D =  \big\{(x,y) \in [0,1]^2 \, : \, x \geq g_2(y) \big\},$ \vspace{-1mm}
\end{center}
are convex.
\end{minipage} \hspace{5mm}
\begin{minipage}{60mm}
\vspace*{4mm}
\includegraphics[scale =.24]{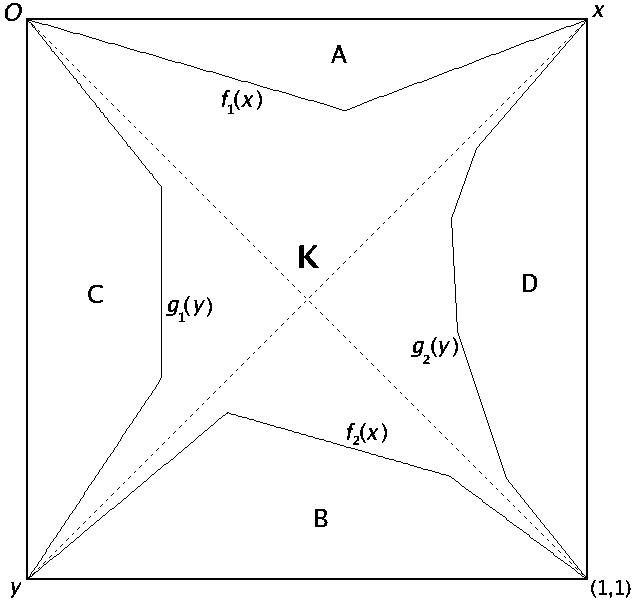} \\
\hspace*{20mm} \vspace{6mm} Fig. 12
\end{minipage}

For short, only a special case is examined; from wich, by extension \\
and symmetry, it follows easily how to treat the general case. \quad Let
$$f_1(x)=\left\{
\begin{array}{lcl}
\dfrac{ax}{b} & \textrm{for} & 0 \leq x \leq \dfrac{bc}{ad+bc}, \vspace{1.5mm} \\
\dfrac{c(1-x)}{d} & \textrm{for} & \dfrac{bc}{ad+bc} \leq x \leq 1,
\end{array}
\right.$$
with integers $\, a \geq 1, \ \, b \geq a, \ \, c \geq 1 \ \, \textrm{and} \ \, d \geq c;$ \\
$f_2(x)=1 \, , \quad g_1(y)=0 \, , \quad \vspace{2mm} g_2(y)=1 \, .$

If on the plane $\, z = 0, \,$ we consider the points \vspace{2mm} $\, P \equiv \Big(\dfrac{bc}{ad+bc}, \dfrac{ac}{ad+bc} \Big) \,$ \\
and $\,Q \equiv (0,1), \,$ then the segments $\, OP \,$ and $\, PQ \,$  represent $\, f_1(x). \,$ \\
In the sheaves $\qquad z-y+\lambda(by-ax)=0 \qquad \textrm{and} \qquad \vspace{2mm} z-y+\mu(dy+cx-c)$ \\
we should choose respectively the \vspace{-1.5mm} planes
$$z=(1-b)y+ax \hspace{35,5mm} (1)$$
$$\textrm{and} \quad z=(1-d)y+c(1-x), \hspace{23mm} (2) \hspace{10mm}$$
wich intersect each other along the straight line
$$ \left\{ \begin{array}{l}
(d-b)y+(a+c)x = c \\
(1-b)y+ax = z.
\end{array} \right. \hspace{34mm}$$
The plane $\, z = x \,$ intersects the plane (1) on the straight line
$$ \left\{ \begin{array}{l}
(b-1)y = (a-1)x \\
z = x.
\end{array} \right. \hspace{26mm}$$
On the plane $\, z = 0, \ $ the straight lines $ \quad (b-1)y = (a-1)x$ \\
and \vspace{3mm} $ \quad (d-b)y+(a+c)x = c \quad $ intersect each other in a point $\, S \equiv (x_{_S}, y_{_S}), \,$ where  \vspace{3mm} $\, x_{_S} = \dfrac{c(b-1)}{(b-1)(a+c)+(d-b)(a-1)}. \,$ \\
Therefore, we must distinguish three cases, according to \vspace{3mm} $\, x_{_S} \;$\raisebox{-1.8mm}{\includegraphics[scale =.8]{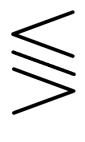}}$\; \dfrac 12.$ \\
If $\, x_{_S} \leq \dfrac 12, \,$ the planes $ \ z = 1-x \ $ and $ \ z = x \ $ intersect the plane (2), respectively on the straight lines
$$ \left\{ \begin{array}{l}
(d-1)y + (c-1)x = c-1 \\
z = 1-x \hspace{37mm} \textrm{and}
\end{array} \right.$$
$$ \left\{ \begin{array}{l}
(d-1)y + (c+1)x = c \\
z = x;
\end{array} \right. \hspace{21mm}$$
so that on the plane $\, z = 0, \,$ the straight line $ \quad (d-1)y + (c+1)x = c \quad $ intersects the straight line \vspace{1mm} $ \quad (d-1)y + (c+1)x = c$ \\
in the point \vspace{1mm} $\, T \equiv \bigg(\dfrac12, \dfrac{1}{2} \cdot \dfrac{c-1}{d-1} \bigg) \,$ \\
and obviously, the straight line $ \quad (b-1)y + (a-1)x = 0 \quad $ in the point $\, S. \,$ \\
So, if $\, R \equiv (0,1/2) \,$ indicates the middle point of $\, OQ, \,$ the triangle $\, OPQ \,$ 
constists of two quadrilaterals $\, OSTR \,$ and $\, PSTQ \,$ and of two triangles $\, OSP \,$ and 
$\, QRT \,$ \quad (Fig. 13).
\begin{center}
\includegraphics[scale =.32]{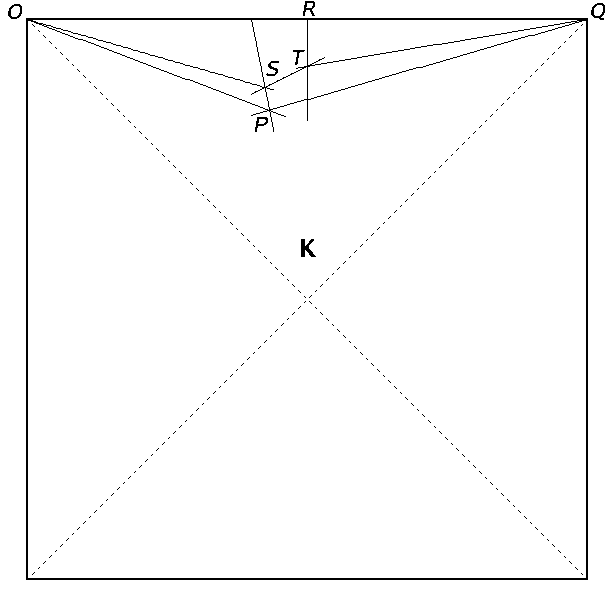} \vspace{1mm} \\
Fig. 13
\end{center}
And again by corollary 2.5 and proposition 3.3, it follows that the functions $\, d_1(x,y) = x$ 
\quad and
$$d_2(x,y)=\left\{
\begin{array}{lll}
(1-b)y+ax & \textrm{if} & (x,y) \in OSP, \\
(1-d)y+c(1-x) \ \; & \, \raisebox{-1.5mm}{''} & (x,y) \in PSTQ, \\
x & \, \raisebox{-1.5mm}{''} & (x,y) \in OSTR, \\
1-x & \, \raisebox{-1.5mm}{''} & (x,y) \in QRT, \\
y & \, \raisebox{-1.5mm}{''} & (x,y) \in \textrm{\textbf{K}}
\end{array}
\right. $$
generate a projective subalgebra of \vspace{3mm} $\, Free_2.$ \\
If $\, x_{_S} = \dfrac 12, \,$ then the points $\, S \,$ and $\, T \,$ coincide; so that, for the 
definition \\
of $\, d_2(x,y), \,$ the four triangles $\, OSP \,$ $\, PSQ \,$ $\, OSR \,$ and $\, QRS \,$ are to 
be considered \quad \vspace{4mm}(Fig. 14).
\begin{center}
\includegraphics[scale =.38]{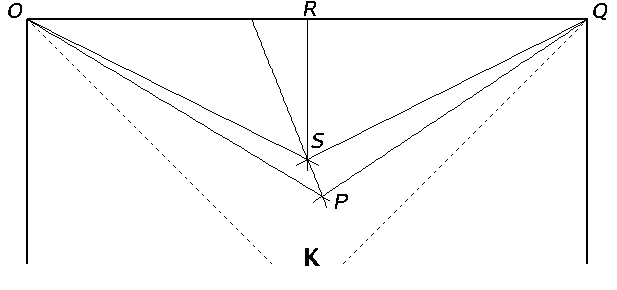} \vspace{1mm} \\
Fig. 14
\end{center}
That is
$$d_2(x,y)=\left\{
\begin{array}{lll}
(1-b)y+ax & \textrm{if} & (x,y) \in OSP, \\
(1-d)y+c(1-x) \ \; & \, \raisebox{-1.5mm}{''} & (x,y) \in PSQ, \\
x & \, \raisebox{-1.5mm}{''} & (x,y) \in OSR, \\
1-x & \, \raisebox{-1.5mm}{''} & (x,y) \in QRS, \\
y & \, \raisebox{-1.5mm}{''} & (x,y) \in \textrm{\textbf{K}}.
\end{array}
\right. \vspace{3mm} $$
Finally, if $\, x_{_S} \geq \dfrac 12, \,$ then on the plane $\, z = 0, \,$ the straight \vspace{2mm} lines \\
$\quad (d-1)y + (c-1)x = c-1 \quad$ and $ \quad (d-b)y+(a+c)x = c \quad $ intersect \vspace{2mm} each other in a point $\ U \equiv (x_{_U}, y_{_U}) \ $ where \vspace{2mm} $\quad x_{_U} = \dfrac{c(b-1)+d-b}{(d-1)(a+c)-(d-b)(c-1)};$ \\
it is easy to check that $\ x_{_S} \geq \dfrac 12 \ $ implies \vspace{2mm} $\ x_{_U} \geq \dfrac 12.$ \\
On the plane $\; z = 0, \;$ the straight line $\quad (b-1)y-(a+1)x = -1, \quad$ projection of the straight line intersection of the plane (1) with the \vspace{1mm} plane $ z = 1-x, \ $ intersects the straight line $ \quad (b-1)y-(a-1)x = 0 \quad $ at the \vspace{1mm} point $\ V \equiv \bigg(\dfrac12, \dfrac{1}{2} \cdot \dfrac{a-1}{b-1} \bigg) \ $ and the straight \vspace{.5mm} line $ \quad (d-1)y + (c-1)x = c-1$ \\
at the point \vspace{1.5mm} $\, U \,$ \quad (Fig. 15). \\
We have:
$$d_2(x,y)=\left\{
\begin{array}{lll}
(1-b)y+ax & \textrm{if} & (x,y) \in OVUP, \\
(1-d)y+c(1-x) \ \; & \, \raisebox{-1.5mm}{''} & (x,y) \in PUQ, \\
x & \, \raisebox{-1.5mm}{''} & (x,y) \in OVR, \\
1-x & \, \raisebox{-1.5mm}{''} & (x,y) \in QUVR, \\
y & \, \raisebox{-1.5mm}{''} & (x,y) \in \textrm{\textbf{K}}.
\end{array}
\right. $$
\newpage \vspace*{4mm} \noindent
\textsc{Example 3.6.} - Let
$$f_1(x)=\left\{
\begin{array}{lcl}
\dfrac{2x}{7} & \textrm{for} & 0 \leq x \leq \dfrac{21}{37}, \vspace{1.5mm} \\
\dfrac{3(1-x)}{8} & \textrm{for} & \dfrac{21}{37} \leq x \leq 1,
\end{array}
\right. \vspace{4mm} $$
$$f_2(x)=1, \qquad g_1(y)=0, \qquad \vspace{6mm} g_2(y)=1.$$
We have: $ \ P \equiv \Big(\dfrac{21}{37},\dfrac{6}{37}\Big), \quad x_{_S} = \dfrac{18}{31}, \quad 
x_{_U} = \dfrac{19}{33}, \quad \vspace{12mm} V \equiv \Big(\dfrac{1}{2},\dfrac{1}{12}\Big); \ $
\begin{center}
\includegraphics[scale =.5]{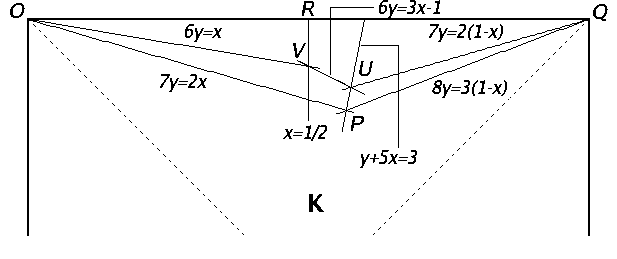} \vspace{1mm} \\
Fig. \vspace{6mm} 15
\end{center}
\hspace*{1mm} $d_1(x,y) = \vspace{2mm} x,$
$$d_2(x,y)=\left\{
\begin{array}{lll}
2(x-3y) & \textrm{if} & y<=2x/7 \ \; \textrm{\textbf{and}} \ \; y<=3-5x \\
& & \textrm{\textbf{and}} \ \; y>=x/6 \ \; \textrm{\textbf{and}} \ \; y>=(3x-1)/6, \\
3(1-x)-7y \ \; & \, \raisebox{-1.5mm}{''} & y>=3-5x \ \; \textrm{\textbf{and}} \ \; y<=3(1-x)/8 \\
& & \textrm{\textbf{and}} \ \; y>=2(1-x)/7, \\
x & \, \raisebox{-1.5mm}{''} & y<=x/6 \ \; \textrm{\textbf{and}} \ \; x<=1/2, \\
1-x & \, \raisebox{-1.5mm}{''} & x>=1/2 \ \; \textrm{\textbf{and}} \ \; y<=(3x-1)/6 \\
& & \textrm{\textbf{and}} \ \; y<=2(1-x)/7, \\
y & \, \raisebox{-1.5mm}{''} & y>=2x/7 \ \; \textrm{\textbf{or}} \ \; y>=3(1-x)/8.
\end{array}
\right. $$
\newpage \noindent
\begin{center}
\includegraphics[scale =.4]{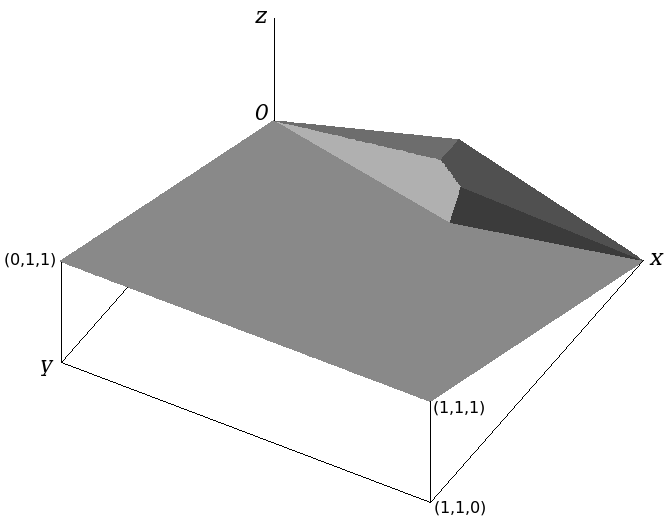} \vspace{-6mm}\\
\small
\hspace*{-40mm} \vspace{3mm} $d_2(x,y)$ \\
\normalsize
Fig. \vspace{6mm} 16
\end{center}
$iii)$ - Now we examine the case that equalizer $\, K \,$ is identified by a generic triangle with a vertex at the origin.
\begin{center}
\includegraphics[scale =.7]{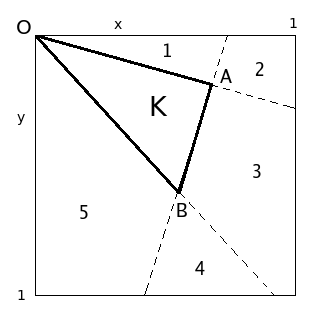} \vspace{-1mm} \\
Fig. 17
\end{center}
Let $\ ax+by = 0, \ \ a_1x+b_1y = 0 \ $ and $\ a_2x+b_2y+c = 0 \ $ be respectively \\
the equations of the straight lines for $\ OA, \ \ OB, \ \ AB \ $ \quad (Fig. 17), \\
with $\ a, \, a_1, \, a_2, \, b, \, b_1, \, b_2 \in \mathbb{Z}, \ \ a, \, a_1, \, c > 0, \ \ b, \, b_1 < 0, \quad \textrm{and} \\
(a,b)=(a_1,b_1)=(a_2,b_2,c)=1.$ \hspace{25mm} Then:
$$K=\{(x,y):ax+by\leq 0 \ \; \textrm{and} \ \; a_1x+b_1y\geq 0 \ \; \textrm{and} \ \; a_2x+b_2y+c\geq 0\}.$$
We can apply theorem 2.4, taking into account the proposition 3.3; \quad more precisely, we 
proceded to the construction of the generators $\, d_1 \,$ and $\, d_2 \,$ assuming that they 
coincide, respectively, with $\, g_1 = x \,$ and $\, g_2 = y \,$ on the equalizer $\, K. \,$ Then, 
we connect these ``basic pieces'' of $\, d_1 \,$ and $\, d_2 \,$ with the plane $\, z = 0, \,$ by 
appropriate planes, so that the hypothesis of the theorem 2.4 are satisfied. \\
For the side $\, OA, \,$ we consider two planes belonging respectively to the sheaves \hspace*{2mm} $\, z=x+k(ax+by) \;$ and $\; z=y+h(ax+by). \,$ \ As the only point of $\, K, \,$ where one of the coordinates $\, x, \, y \,$ takes zero value is the origin of axes, then these two planes must intersect necessary the plane $\, z=0 \,$ on the same straight line. This implies
$$(ka+1)(hb+1)-hakb=0,$$
that is
$$ka+hb+1=0.$$
So, the planes of the relative sheaves that satisfy such condition are:
$$ \left.
\begin{array}{ll}
z=-hbx+kby &  \hspace{20mm} (1) \\
z= \hspace*{2.7mm} hax-kay & \hspace{20mm} (2)
\end{array}
\right.$$
and their intersection straight line with the plane $\, z = 0 \,$ is
$$y=\dfrac hkx. \hspace{38mm} (3)$$
If $\, k_0, \, h_0 \,$ is a particular solution of the diophantine equation  $\, ak+bh=-1, \,$
the general solution has the form:
$$k=k_0 - sb, \ \ h=h_0+sa, \ \ \ \textrm{with} \ \; s\in \mathbb{Z} \ \; \textrm{wichever};$$
so that $\ \ \dfrac{h}{k} = - \dfrac{a}{b} + \dfrac{1/b^2}{s-k_0/b} \ \ \ $ and for $\; s < \dfrac{k_0}{b} + \dfrac{1}{ab} \;$ we have: \vspace{4mm} $\ \ \ 0 < \dfrac hk < - \dfrac ab.$ \\
In a similar manner, we proceed for the side $\, OB; \,$ fixing the planes
$$ \left.
\begin{array}{ll}
z=-h'b_1x+k'b_1y &  \hspace{20mm} (4) \\
z= \hspace*{2.7mm} h'a_1x-k'a_1y & \hspace{20mm} (5)
\end{array}
\right.$$
as well as their intersection straight line with the plane $\, z = 0$
$$y=\dfrac{h'}{k'}x \hspace{44mm} (6)$$
with the conditions $\; \  k'a_1+h'b_1+1=0 \ \;$ and \vspace{3mm} $\; \ \dfrac{h'}{k'} > - \dfrac{a_1}{b_1}.$ \\
Now, we indicate with $\, A' \,$ and $\, A'' \,$, respectively, the intersection points \\
of the vertical for $\, A \,$ with the planes $\, z = x \,$ and $\, z = y \, ; \,$ with $\, B' \,$ and $\, B'' \,$ the intersection points of the \vspace{2mm} vertical for $\, B \,$ with the same planes. \\
Set \ $\delta \, = \, ba_2-ab_2$ \ and \ $\delta_1 \, = \, b_1a_2-a_1b_2$, \ we have:
$$A \ \equiv \ \Big( \dfrac{-bc}{\delta} \ , \, \dfrac{ac}{\delta} \, \Big) \qquad \textrm{and} \qquad B \ \equiv \ \Big( \dfrac{-b_1c}{\delta_1} \ , \, \dfrac{a_1c}{\delta_1} \, \Big),$$
with $\ \delta, \ \delta_1 \ > 0; \ \ \ $ consequently
$$A' \ \equiv \ \Big( \dfrac{-bc}{\delta} \ , \, \dfrac{ac}{\delta} \ , \, \dfrac{-bc}{\delta} \,\Big), \ \qquad A'' \ \equiv \ \Big( \dfrac{-bc}{\delta} \ , \, \dfrac{ac}{\delta} \ , \, \dfrac{ac}{\delta} \,\Big) \qquad \textrm{and}$$
$$B' \ \equiv \ \Big( \dfrac{-b_1c}{\delta_1} \ , \, \dfrac{a_1c}{\delta_1} \ , \, \dfrac{-b_1c}{\delta_1} \,\Big), \qquad B'' \ \equiv \ \Big( \dfrac{-b_1c}{\delta_1} \ , \, \dfrac{a_1c}{\delta_1} \ , \, \dfrac{a_1c}{\delta_1} \,\Big). \vspace{5mm} \quad$$
We consider the sheaves of \vspace{-2mm} planes
$$z-x+l(a_2x+b_2y+c)=0 \qquad \textrm{and} \qquad z-y+m(a_2x+b_2y+c)=0,$$
generated by the plane $a_2x+b_2y+c=0 \,$ and, respectively, by the planes $\, z=x$ and $\, z=y; \,$ \ these sheaves intercept on the plane $\,z=0 \,$ the sheaves of straight lines
$$x=l(a_2x+b_2y+c) \qquad \textrm{and} \qquad y=m(a_2x+b_2y+c).$$
So the coordinates of the point $\, P, \,$ intersection of a straight line of the first sheaf with one straight line of the second, must satisfy the system
\begin{displaymath}
\left\{ \begin{array}{rl}
l(a_2x+b_2y+c) & = \ x \\
m(a_2x+b_2y+c) & = \ y.
\end{array} \right.
\end{displaymath}
It follows, immediately, that the point $\, P \,$ must lie on the straight line through the origin of equation \ $mx=ly$ \ and that
$$ P\equiv \Big(\dfrac{lc}{1-la_2-mb_2} \ , \, \dfrac{mc}{1-la_2-mb_2}\Big).$$
In order that, the straight line through $\, OP \,$ lies between the straight \vspace{1.5mm} lines 
through $\, OA \,$ and $\, OB \,$ \quad (Fig. 18), \quad it must be $\ -\dfrac{a}{b}<\dfrac{m}{l}<-\dfrac{a_1}{b_1},$ \vspace{.5mm} \\
and the parameters $\, m \,$ and $\, l \,$ must have the same sign. \\
If $\ m < 0 \ $ and $\ l < 0 \ $, \ must be also $\ 1-la_2-mb_2 < 0; \ $ so that is
$$a_2\dfrac{lc}{1-la_2-mb_2}+b_2\dfrac{mc}{1-la_2-mb_2}+c=\dfrac{c}{1-la_2-mb_2}<0,$$
therefore  $\, P \,$ is an outer point to the equalizer $\, K \,$ (it is in zone 3) and when $\, m \,$ and $\, l \,$ in absolute value increase, its distance from the side $\, AB \,$ of $\, K \,$ decreases. \\
The equation of a generic straight line of the sheaf generated by a straight line of the first sheaf with one of the second, (that is a straight line of the plane $\, z = 0 \,$ passing through $\, P) \,$ is:
$$h\big[l(a_2x+b_2y+c)-x\big]+k\big[m(a_2x+b_2y+c)-y\big]=0,$$
with $\ h, \ k \, \in \mathbb{Z} \ $, wich can be put in the form
$$\big[a_2(hl+km)-h\big]x+\big[b_2(hl+km)-k\big]y+c(hl+km)=0.$$
Let
$$ \left.
\begin{array}{lll}
& \bar{a}x+\bar{b}y+\bar{c} = 0 &  \hspace{20mm} (7) \\
\textrm{and} \qquad & \hat{a}x+\hat{b}y+\hat{c} = 0 & \hspace{20mm} (8)
\end{array}
\right.$$
be the equations for two straight lines belonging to this sheaf. \\
The equation of a McNaughton generic plane of the sheaf whose axis is the straight line (7), can be written in the form
$$z=\lambda(\bar{a}x+\bar{b}y+\bar{c})=\lambda\Big\{\bar{h}\big[l(a_2x+b_2y+c)-x\big]+\bar{k}\big[m(a_2x+b_2y+c)-y\big]\Big\},$$
with $\lambda$ integer, \ and in order that the point $\, A' \,$ lies on this plane, it must be:
$$-\dfrac{bc}{\delta}=\lambda' \Big\{\bar{h}\Big[\dfrac{bc}{\delta}\Big]+\bar{k}\Big[ -\dfrac{ac}{\delta} \Big]\Big\}, \quad \textrm{from \ wich} \quad -b=\lambda' (\bar{h}b -\bar{k}a);$$
while, if we require that the point $\, A'' \,$ lies on a plane of the same sheaf, \\
it must be:
$$\dfrac{ac}{\delta}=\lambda'' \Big\{\bar{h}\Big[\dfrac{bc}{\delta}\Big]+\bar{k}\Big[-\dfrac{ac}{\delta} \Big]\Big\}, \quad \textrm{from \ wich} \quad a=\lambda'' (\bar{h}b-\bar{k}a).$$
As $\ (a,b)=1, \ $ it follows \vspace{-1.5mm} necessarily
$$\bar{h}b-\bar{k}a=\mp 1, \qquad \lambda'=\pm b \qquad \textrm{and} \qquad \lambda'' =\mp a.$$
If we choose $\, \bar{h}b-\bar{k}a = 1, \,$ then $\, \lambda' = -b \,$ and $\, \lambda'' = a \,$ and the equations of two planes of the sheaf, one passing through $\, A' \! , \,$ the other through $\, A'' \! \! , \,$ are respectively:
$$ \left.
\begin{array}{lll}
& z = -b\bar{a}x-b\bar{b}y-b\bar{c} &  \hspace{20mm} (9) \\
\textrm{and} \qquad & z = \hspace*{2.7mm} a\bar{a}x+a \bar{b}y+a\bar{c}. & \hspace{18.2mm} (10)
\end{array}
\right.$$
These planes intersect $\, z$-axis in two points with \\
\hspace*{12mm} $\, z' = -b\bar{c} = -bc(\bar{h}l+\bar{k}m) \qquad \textrm{and} \qquad z'' = a\bar{c} = ac(\bar{h}l+\bar{k}m);$ \\
so that $\, z' \ \textrm{and} \ z'' > 0 \,$ if and only if $\, \bar{h}l+\bar{k}m > 0.$ \\
Now the equation $\ \bar{h}b-\bar{k}a = 1 \ $ admits the integer solutions: \\
$\bar{h} = h_1+at, \ \bar{k} = k_1+bt, \ $ with arbitrary integer parameter $\, t; \,$ so the condition $\, \bar{h}l + \bar{k}m > 0 \,$ is equivalent to $\, t(la+mb) > -mk_1-lh_1 \,$ and this is \vspace{2.5mm} satisfied for $\, t > - \dfrac{mk_1+lh_1}{la+mb} \,$ (being \vspace{3.5mm} $\, la+mb > 0).$ \\
Similarly, the equations of the two planes of the sheaf with the straight line (8) as axis, passing respectively through the points $\, B' \,$ and $\, B'', \,$ are:
$$ \left.
\begin{array}{lll}
& z = -b_1\hat{a}x-b_1\hat{b}y-b_1\hat{c} &  \hspace{20mm} (11) \\
\textrm{and} \qquad & z = \hspace*{2.7mm} a_1\hat{a}x+a_1 \hat{b}y+a_1\hat{c} \; ; & \hspace{20mm} (12)
\end{array}
\right.$$
with the integer parameters $\, \hat{h} \,$ and $\,\hat{k} \,$ (wich determine $\, \hat{a}, \ \; \hat{b} \,$ and $\,\hat{c}) \,$ satisfying \\
the conditions $\ \hat{h}b_1-\hat{k}a_1 = 1, \ $ and \vspace{2.5mm} $\  \hat{h}l+\hat{k}m > 0.$ \\
Finally, we note that the planes through the points $\ A'\!, \ \; B'\!\!, \ \; P \ \ \textrm{and} \\
A''\!\!, \ \; B''\!\!\!, \ \; P \ $ have, respectively, \vspace{-1.5mm} equations:
$$ \left.
\begin{array}{lll}
& z = (1-la_2)\: x-lb_2\: y-lc &  \hspace{20mm} (13) \\
\textrm{and} \qquad & z = -ma_2\: x+(1-mb_2)\: y-mc. & \hspace{20mm} (14)
\end{array}
\right.$$
We denote, also respectively, by $\, Q \,$ and $\, R \,$ the intersection points between the straight lines (3), (7) and (6), (8).
\begin{center}
\includegraphics[scale =.9]{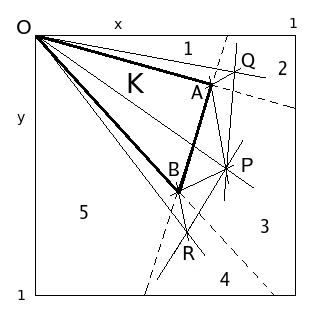} \vspace{-5mm} \\
Fig. 18
\end{center}
So the generators $\, d_1 \,$ and $\, d_2, \,$ \\
\underline{zero} along and outside of the broken line $\, O \; Q \; P \; R \; O, \quad$ coincide: \\
\hspace*{10mm} $\ d_1 \,$ with planes\, (1),\, (9),\, (13),\, (11),\, (4)\, and $\, z = x$ \\
and $\quad d_2 \,$ with planes (2),\, (10),\, (14),\, (12),\, (5)\, and $\, z = y,$ \\
respectively, on the triangles \\
$\, O \; Q \; A, \ \ Q \; A \; P, \ \ A \; P \; B, \ \ P \; B \; R, \ \ B \; R \; O \;$ and 
\textit{\textbf{K}} \quad (Fig. \vspace{1mm} 18). \\
It follows immediately, by proposition 3.3, that the generators $\, d_1 \,$ and $\, d_2, \,$ \\
so defined, satisfy the conditions of the theorem \vspace{1.5mm} 2.4. \\
These considerations can be easily extended also to the case that \\
equalizer $\, K \,$ is decomposable into triangles, all with a vertex in the origin. Infact, if the construction, above described, is repeated for each triangle
$$\, O \; A_1 \; B_1, \; \ O \; A_2 \; B_2, \; \dots, \; \ O \; A_n \; B_n, \,$$
we obtain in particular respectively the quadrilaterals
$$\, O \; Q_1 \; P_1 \; R_1, \; \ O \; Q_2 \; P_2 \; R_2, \; \dots, \; \ O \; Q_n \; P_n \; R_n, \,$$
whose union is a polygon
$$\, O \; Q_1 \; P_1 \; S_1 \; \dots \; S_m \; P_n \; R_n \; O \,$$
that containes $\, K. \,$ \ As each segment of the perimeter of this polygon belongs to one of the 
quadrilaterals $\, O \; Q_i \; P_i \; R_i, \,$ it is possible to use convenient pieces of planes 
determinated for each triangle, to define the \vspace{4mm} generators (Fig. 19).
\begin{center}
\includegraphics[scale =.4]{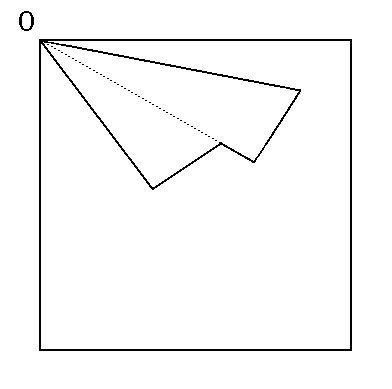}
\includegraphics[scale =.4]{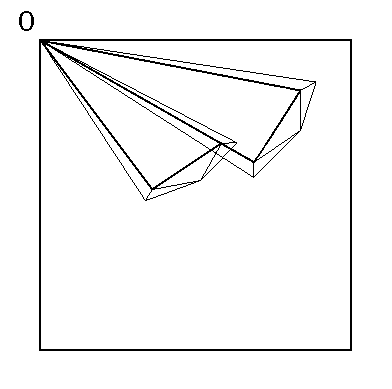} \vspace{-5mm} \\
\includegraphics[scale =.4]{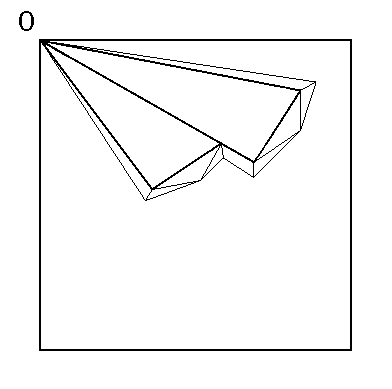} \vspace{-5mm} \\
Fig. 19 \\
\end{center}
In conclusion, using the methods described in \S\,3 $b),$ in the case of equalizers containing the diagonals and whose boundary consists of continuous and piecewise linear functions, projective generators are obtained. Similarly, the same result is obtained in \S\,3 $c),$ in the case of equalizers formed by triangles with a common vertex coinciding with one of the vertices of $\, [0,1]^2.$ Apart from the complication in the calculations, the results can be extended to the case of three or more generators. It also seems likely that the methods described in \S\,3 can be applied in other \vspace{50mm} cases. \\
\textsc{References} \vspace{2mm}
\small
\setlength{\leftmargini}{1.25em}
\begin{enumerate}
\item
\textsc{L. M. Cabrer - D. Mundici}, \textit{Projective MV-algebras and rational polyhedra}, Algebra Universalis \textbf{62}, no. 1, 63-74 (2009).
\item
\textsc{C. C. Chang}, \textit{Algebraic analysis of many valued logics}, Trans. Amer. Math. Soc., \textbf{88}, 467-490 (1958).
\item
\textsc{C. C. Chang}, \textit{A new proof of the completeness of the {\L}ukasiewicz axioms}, Trans. Amer. Math. Soc., \textbf{93}, 74-80 (1959).
\item
\textsc{R. Cignoli - I. D'Ottaviano - D. Mundici},  Algebraic foundations of many-valued reasoning, Trends in Logic, Kluwer, Dordrecht 2000.
\item
\textsc{A. Di Nola - R. Grigolia - A. Lettieri}, \textit{Projective MV-algebras}, Internat. J. Approx. Reason. \textbf{47}, no. 3, (2008) 323-332.
\item
\textsc{A. Di Nola - R. Grigolia}, \textit{Projective MV-algebras and their automorphism groups}, J. Mult.-Valued Logic Soft Comput. \textbf{9}, no. 3, 291-317 (2003).
\item
\textsc{F. Lacava}, \textit{A characterization of free generating sets in MV-algebras}, Algebra Universalis \textbf{57}, no. 4, 455-462 (2007).
\end{enumerate}
\end{document}